\documentclass{article} 
\usepackage{amsfonts,amsmath,amssymb,amstext,amsbsy,amsopn}
\usepackage{url}
\usepackage{graphicx}
\usepackage{booktabs}
\usepackage{subfigure}
\date{}

\newcommand{\R}{{\mathbb R}}
\newcommand{\Z}{{\mathbb Z}}

\newcommand{\N}{{\mathbf N}}
\newcommand{\itg}{\int \limits}
\newcommand\ds{\displaystyle}
\newcommand{\de}{\partial}
\newcommand{\ee}{{\rm e}}
\newcommand{\supp}{{\rm supp\,}}
\newcommand{\im}{{\rm Im\,}}
\newcommand{\re}{{\rm Re\,}}

\newcommand{\erfc}{{\rm erfc}}

\def\a{{\alpha}}

\def\eps{\epsilon}
\def\e{\varepsilon}

\def\D{\Delta}
\def\t{\tau}
\def\h{\eta}

\def\x{\xi}

\def\O{\Omega}
\def\o{\omega}

\def\k{\kappa}

\newcommand{\ba}{\boldsymbol{\alpha}}

\newcommand{\cD}{\mathcal D}

\newcommand{\cK}{\mathcal K}

\newcommand{\cM}{\mathcal M}
\newcommand{\cN}{\mathcal N}

\newcommand{\cO}{\mathcal O}
\newcommand{\cP}{\mathcal P}
\newcommand{\cQ}{\mathcal Q}
\newcommand{\cR}{\mathcal R}
\newcommand{\cS}{\mathcal S}

\newcommand{\cF}{\mathcal F}

\newcommand{\bx}{{\mathbf x}}
\newcommand{\by}{{\mathbf y}}

\newcommand{\bm}{{\mathbf m}}
\newcommand{\bk}{{\mathbf k}}
\newcommand{\bP}{{\mathbf P}}
\newcommand{\bQ}{{\mathbf Q}}

\newtheorem{thm}{Theorem}[section]
\numberwithin{thm}{section}
\newtheorem{remark}{Remark}[section]
\numberwithin{remark}{section}

\title{A fast solution method for time dependent multidimensional Schr\"odinger equations}
\author{{  F. Lanzara$^{\mbox{\tiny 1}}$ , V. Maz'ya$^{\mbox{\tiny 2}}$ ,
G. Schmidt$^{\mbox{\tiny 3}}$}}

\begin{document}
\maketitle

\parbox{10cm}{\begin{flushleft}
{\footnotesize\em
\begin{itemize}
\item[$^{\mbox{\tiny\rm 1}}$] Department of Mathematics, Sapienza University of Rome,\\
Piazzale Aldo Moro 2, 00185 Rome, Italy\\
\texttt{\rm lanzara\symbol{'100}mat.uniroma1.it}
\item[$^{\mbox{\tiny\rm 2}}$]Department of Mathematics, University of
Link\"oping, \\ 581 83 Link\"oping, Sweden\\
\texttt{\rm vlmaz\symbol{'100}mai.liu.se }
\item[$^{\mbox{\tiny\rm 3}}$]Weierstrass Institute for Applied Analysis and
Stochastics, \\  Mohrenstr. 39,
10117 Berlin, Germany \\
\texttt{\rm schmidt\symbol{'100}wias-berlin.de}
\end{itemize}
}
\end{flushleft}}

{\bf{Abstract:}} In this paper we propose fast solution
methods for the Cauchy problem for the multidimensional Schr\"odinger
equation. Our approach is based on the approximation of the data by
the basis functions introduced in the theory of approximate
approximations.  We obtain high-order approximations also in higher
dimensions up to a small saturation error, which is  negligible
in computations, and we prove error estimates in mixed Lebesgue spaces
for the inhomogeneous equation. The proposed method is very efficient
in high dimensions if the densities allow separated
representations. We illustrate the efficiency of the procedure on
different examples, up to approximation order $6$ and space dimension
$200$.



\section{Introduction}\label{sec1}
\setcounter{equation}{0}

The present paper is devoted to the numerical solution of initial value problems for
the Schr\"odinger equation of  free
particles
\begin{equation}\label{main0} i \frac{\partial u}{\partial t} +\D_\bx
u = 0 \,,
\end{equation} where $u=u(\bx,t)$ is the wave function depending on
the spatial variables $\bx=(x_1,...,x_n)\in\R^n$ and the time
$t\in\R$, $\D_\bx$ is the usual Laplacian with respect to the
variables $\bx$. The time evolution of physical systems is generally
described via partial differential equations, especially via the
Schr\"odinger equation 
\[
i  \frac{\partial u}{\partial t}+ \D_\bx u= V(\bx,t) u\,,
\]
in such fields where wave propagation is considered, for example,
optics, acoustics and quantum mechanics. Here $V (\bx,t)$ is the
potential, which models the interaction of the particle with its
environment.  In case of a free particle $V(\bx,t)=0$ and we get
\eqref{main0}.  The Schr\"odinger equation with real valued time
dependent potential $V(t)$
\[
i\, \frac{\de w}{\de t}+\D_\bx w=V(t) w
\]
can be dealt with \eqref{main0} by the transformation $w(\bx,t)=u(\bx,t) \ee^{-i\int_0^t V(\t)d\t}$.

Solving numerically the Schr\"odinger equation is of great practical
use but is in general quite a complex problem mainly due to the fact
that the wave function propagates high frequency oscillations.  
The main objective of the paper  is to develop  fast solution methods for
the Cauchy problem of \eqref{main0} 
\begin{equation}\label{initcond}
u(\bx,0)=g(\bx),\quad \bx=(x_1,,\ldots,x_n)\in\R^n,
\end{equation}
and for the inhomogeneous Schr\"odinger equation
\begin{equation}\label{inhom}
i \frac{\partial u}{\partial t} +\D_\bx u=f(\bx,t),\quad
(\bx,t)\in \R^n\times \R_+ \, 
\end{equation}
which is effective also in high dimension.

Under suitable integrability or decay conditions on $g$ and $f$
the solution of \eqref{main0} - \eqref{initcond} can be written as
\begin{equation}\label{sol0}
u(\bx,t) = \cS g(\bx,t)=  
\itg_{\R^n} \cK(\bx-\by,t)\,g(\by)\, d\by \, ,
\end{equation}
and the solution of \eqref{inhom} - \eqref{initcond} is given by
\begin{equation*}\label{solinh}
u(\bx,t)= \cS g (\bx,t)+\Pi  f(\bx,t)
\end{equation*}
with
\begin{equation}\label{sol1}
\Pi  f(\bx,t)=-i \itg_0^{t} ds \itg_{\R^n} \cK(\bx-\by,t-s) f(\by,s)\, d\by \,ds 
=-i\itg_0^t (\cS f(\cdot,s))(\bx,t-s)\, ds \, .
\end{equation}
Here $\cK(\bx,t)$ denotes  the fundamental solution of \eqref{main0} \cite[p.193]{LCE}
\begin{equation*} \label{fundsol}
\cK(\bx,t)=\frac{\ee^{ \, i \, |\bx|^2/(4t)}}{(4 \pi i t)^{n/2}}  \, .
\end{equation*}
In terms of the Fourier transform 
\begin{equation*} \label{four}
\cF u (\boldsymbol{\xi}) = \itg_{\R^n} u(\bx) \, 
\ee^{-2 \pi i \langle \bx, {\boldsymbol{\xi}} \rangle } \, d \bx 
\end{equation*}
the integral operators $\cS$ and $\Pi$ can be written as follows:
\begin{align}\label{solF}
\cS g(\bx,t) &=  \itg_{\R^n}\ee^{\,2 \pi i \langle \bx, \boldsymbol{\xi}\rangle} 
\ee^{-4 \pi^2 i t|\boldsymbol{\xi}|^2} \cF g(\boldsymbol{\xi})\,d\boldsymbol{\xi} \, ,\\\nonumber
\Pi  f(\bx,t)&= -i \itg_0^{t} ds \itg_{\R^n} \ee^{\,2 \pi i \langle \bx, \boldsymbol{\xi}\rangle} 
\ee^{-4 \pi^2 i \, (t-s)|\boldsymbol{\xi}|^2} \cF f(\boldsymbol{\xi},s)\,d\boldsymbol{\xi}\, .
\end{align}
For fixed $t>0$ the integral operator $\cS$
is bounded in $L^p=L^p(\R^n)$ spaces. We write $||f||_{L^p}$ the Lebesgue norm of a function $f\in L^p$.
From  \eqref{sol0} and \eqref{solF} it follows immediately that
\begin{align*}
\|\cS g(\cdot,t)\|_{L^\infty}\leq \frac{1}{(4 \pi |t|)^{n/2}} \|g\|_{L^1} \, ,
\quad \|\cS g(\cdot,t)\|_{L^2} =  \|g\|_{L^2} \, ,
\end{align*}
hence by interpolation 
the  $L^p$ dispersive estimate  %{\color{blue} \cite{CZ}}
\begin{align}\label{estip}
\|\cS g(\cdot,t)\|_{L^p} \le 
C |t|^{-n (1/2-1/p)}\|g\|_{L^{p'}} \, ,
\quad t\neq 0\, ,
\end{align}
holds for $2 \le p \le \infty$ and $p'$ is the adjoint exponent, $1/p+1/p'=1$.

Estimates of norms of solutions $u(\bx,t)$ on $\R^n \times \R$
are known for example 
in mixed Lebesgue spaces.
For an interval $I$  and $r,q \ge 1$,
$L^{r,q}(I)$ denotes the Banach space of $L^r(\R^n)$-valued
$q$-summable functions over  $I$ 
with the norm
   \begin{equation*} \label{defmixnorm}
\|u\|_{L^{r,q}(I)}=\|u\|_{r,q}= 
\bigg(\itg_{I} \Big( \itg_{\R^n} |u(\bx,t)|^{r} d\bx\Big)^{q/r} dt \bigg)^{1/q} \, .
   \end{equation*}
The exponent pair $(q,r)$ is called Schr\"odinger-admissible if  
$q,r \ge 2$, $(q,r) \ne (2,\infty)$ and
   \begin{align} \label{defadmissible}
\frac{2}{q}+  \frac{n}{r} = \frac{n}{2} \, .
\end{align}
For any Schr\"odinger-admissible pairs $(q,r)$  Strichartz type estimates
   \begin{align} \label{striest}
\| \cS g\|_{L^{r,q}(\R_+)}  \le C \| g\|_{L^2} \, , \quad
\| \Pi  f\|_{L^{r,q}(\R_+)}\le C 
\|f\|_{L^{r',q'}(\R_+)}
   \end{align}
are valid with constants $C$ independent of $g \in L^2(\R^n)$ and $f \in  L^{r',q'}(\R_+)$, \cite[(11) ]{Tao}. 
Moreover, $u= \cS g+\Pi  f$ is continuous in $t$ in the space $L^2$ and
   \begin{align*} 
\sup_{t \in \R_+} \|u(\cdot ,t)\|_{L^2}  \le C (\| g\|_{L^2} +\|f\|_{L^{r',q'}(\R_+)}) \, . 
   \end{align*}
Note that for $q=r=2+4/n$
we derive the classical Strichartz estimate \cite[(3.2)]{Stri}
   \begin{align*} 
\|u \|_{L^{2(n+2)/n}}=\bigg(\itg_{\R_+}\itg_{\R^n} |u(\bx,t)|^{2(n+2)/n} d\bx dt \bigg)^{n/(2(n+2))}\!\!\!\!\!\!\!\! \le C \big(
\|g\|_{L^2} + \|f\|_{L^{2(n+2)/(n+4)}} \big)\, .
   \end{align*}

The goal of this paper is to derive semi-analytic cubature formulas
for $\cS g$ in \eqref{sol0} and $\Pi f$ in \eqref{sol1} of an
arbitrary high-order which are fast and accurate also if the space
dimension $n\geq 3$. We follow the philosophy introduced in
\cite{LMS2} and \cite{LMS3} for the cubature of high-dimensional
Newton potential over the full space and over half-spaces. The idea is
to approximate the density functions by the basis functions introduced
in the theory of {\it approximate approximations} (cf. \cite{MSbook}
and the references therein). This approach, combined with {\it
separated representations} (cf. \cite{BM} and \cite{BM2}) makes the
method fast and successful  in high dimensions. In \cite{LMS4} and \cite{LMS5} we applied
this procedure to obtain cubature formulas for advection-diffusion
operators over rectangular boxes in $\R^n$. In \cite{LMS6} 
 our approach was extended to parabolic problems. For the
Schr\"odinger equation the situation is different because the
fundamental solution does not decay exponentially and standard
cubature methods are very expansive due to the oscillations of the
kernel, especially in multidimensional case. The application of
approximate approximations to this equation 
reduces these problems and
provides new very efficient semi-analytic cubature formulas.

The article is organized as follows. In section \ref{sec2}, after an
introduction into simple cubature formulas for the operators $\cS$ and $\Pi$
based on  {\it approximate quasi-interpolants}, we prove new
estimates of the cubature error for general generating functions in mixed Lebesgue spaces.
Similar to other integral operators of potential theory, we
obtain high-order approximations also in higher dimensions up to a small
saturation error, which is  negligible in computations.
In section \ref{sec3}
we describe algorithms for 
high-order approximations of \eqref{sol0} and \eqref{sol1}.
Using the tensor product structure of cubature formulas,
these algorithms are very efficient in high dimensions,
if $g$ and $f$ allow  separate representations.
The approach is extended in section \ref{sec4}  to the case
that $g$ and $f$
are supported with respect to $\bx$ in a hyper-rectangle on
$\R^n$. In Section \ref{sec5} we illustrate the
efficiency of the method on several examples, up to approximation
order $6$ and space dimension $200$. For the two-dimensional initial
value problem \eqref{main0}-\eqref{initcond} we provide graphics of the evolution of
 $u(\bx,t)$.  %Section \ref{sec6} is dedicated to a brief conclusion.

\section{Cubature of $\cS g$ and $\Pi f$}\label{sec2}
\setcounter{equation}{0}

\subsection{Approximate quasi-interpolants}
To find an approximate solution of \eqref{main0} - \eqref{initcond}
we replace the function $g$  in \eqref{sol0} by an approximate quasi-interpolant
\begin{equation}\label{appN}
(\cM_{h\sqrt{\cD}} \,g) (\bx)=\cD^{-n/2} \sum_{\bm \in \Z^n} g(h \bm) \h\Big(\frac{\bx-h\bm}{h\sqrt{\cD}}\Big),
\end{equation}
where $\h$ is a rapidly decaying function of the Schwarz space $S(\R^n)$
satisfying for positive integer $N$ the moment condition
   \begin{equation}  \label{momcon} 
    \itg _{{\R} ^n} \eta(\bx) \, d \bx = 1
  \; , \;                 
    \itg _{{\R} ^n} 
     {\bx}^{\boldsymbol{\a}} \eta(\bx) \, d \bx = 0 , \quad
    \forall \, \boldsymbol{\a}  \; , \; 
   1  \le |\boldsymbol{\a}| <  N \, .
    \end{equation}
Here and in the following we
use  multi-index notation, bold Greek letters denote multi-indices.
Then
the function
\begin{equation}\label{appsol0}
\begin{split}
\cS_h g(\bx,t) &= \cS(\cM_{h\sqrt{\cD}} \,g) (\bx,t)\\
&= \frac{1}{\cD^{n/2} (4 \pi i t)^{n/2}} \sum_{\bm \in \Z^n} g(h \bm)
  \itg_{\R^n}\ee^{\, i \, |\bx-\by|^2/(4t)}\h\Big(\frac{\by-h\bm}{h\sqrt{\cD}}\Big)d\by
\end{split}
\end{equation}
can be considered as cubature of $\cS g$, if $\h$ is chosen such
that $\cS \h(\bx,t)$ can be computed  easily, preferably as
an analytic expression. The existence of those generating functions $\h$
has been shown for various integral operators.

Then the cubature error follows immediately 
from the quasi-interpolation error due to
\begin{equation*} \label{simpS}
\cS g(\cdot,t)- \cS_h g(\cdot,t)=\cS(I-\cM_{h\sqrt{\cD}}) \,g (\cdot,t) \, .
\end{equation*}
Approximation properties of quasi-interpolants of the form \eqref{appN}
have been studied in the framework of approximate approximations
(cf. \cite{MSbook}). 
Let us recall the structure of the quasi-interpolation error, which is  proved
in general form  in  \cite[Thm 2.28]{MSbook}.
Suppose that $g$ has generalized derivatives of order $N$.
Using Taylor expansions of $g(\bx)$ for the nodes $h\bm$, $\bm \in {\Z}^n$,
and Poisson's summation formula the quasi-interpolant can be written as
\begin{align}\label{repre}
(\cM_{h\sqrt{\cD}} \,g) (\bx) = &(-h\sqrt{\cD})^N g_N(\bx)
+\sum _{ |\boldsymbol{\alpha}|  = 0} ^{N-1}
\frac{(h \sqrt{\cD})^{|\boldsymbol{\alpha}|}}{\boldsymbol{\alpha} !(2 \pi i)^{|{\boldsymbol{\alpha}}|}} 
\, \partial ^{\boldsymbol{\alpha}} g(\bx)  \,    \sigma_{\boldsymbol{\alpha}}(\bx,\h,\cD)
\end{align}
    with the function 
\begin{align*} 
g_N(\bx)\!=\! \frac{1 }{{\cD}^{n/2}}\!\!
 \sum _{ |\boldsymbol{\alpha}|  = N} \frac{N}{\ba!} \!
 \sum_{\bm \in {\Z}^n} \!\!
\Big(\frac{ \bx-h \bm }{h\sqrt{\cD}}\Big)^{\boldsymbol{\alpha}}  \!
\h \Big(
\frac{\bx\!-\!h\bm}{h\sqrt{\cD}} \Big)
\!     \itg _0 ^1 \! s^{N-1} \partial ^{\boldsymbol{\alpha}} g
     (s \bx + (1-s) h\bm)  \, ds ,
\end{align*}
containing  the remainder of the Taylor expansions,
and the fast oscillating functions
   \begin{align} \label{defsigma}
     \sigma_{\boldsymbol{\alpha}}(\bx,\h,\cD)= \frac{1 }{{\cD}^{n/2}} \sum_{\bm \in {\Z}^n} 
\Big(\frac{ \bx-h \bm }{h\sqrt{\cD}}\Big)^{\boldsymbol{\alpha}}  \!
\h \Big(
\frac{\bx\!-\!h\bm}{h\sqrt{\cD}} \Big)=\sum_{\boldsymbol{\nu} \in {\Z}^n} 
\partial^{\boldsymbol{\alpha}}\cF \eta(\sqrt{\cD} {\boldsymbol{\nu}})
\, \ee^{ \frac{2 \pi i}{h} \langle {\mathbf x}, {\boldsymbol{\nu}}\rangle } \, .
   \end{align}
If  $g \in W_p^N(\R^n)$ with $N>n/p$, $1 \le p \le \infty$, then $g_N$ can be estimated by
   \begin{equation*}\label{ch1rem}
\|g_N\|_{L_p} \le C_N \sum _{ |\boldsymbol{\alpha}|  = N} \|\partial^{\boldsymbol{\alpha}}g\|_{L_p}=
C_N  |g|_{W_p^N}
   \end{equation*}
with a constant $C_N$ depending only on $\h$, $n$, and $p$.
It follows from \eqref{defsigma} that due to the moment condition \eqref{momcon} the second sum in \eqref{repre}
transforms to
   \begin{align*} 
g(\bx)+\sum _{ |\boldsymbol{\alpha}|  = 0} ^{N-1}
\frac{( h \sqrt{\cD})^{|\boldsymbol{\alpha}|}}{\boldsymbol{\alpha} !(2 \pi i)^{|{\boldsymbol{\alpha}}|}} 
\, \partial ^{\boldsymbol{\alpha}} g(\bx) \, \eps_{\boldsymbol{\alpha}}(\bx,\h,\cD) \, ,
   \end{align*}
where we denote 
   \begin{align*} 
\eps_{\boldsymbol{\alpha}}(\bx,\h,\cD) = \sum_{\substack{\boldsymbol{\nu} \in {\Z}^n\\ \boldsymbol{\nu}\neq 0}} 
\partial^{\boldsymbol{\alpha}}\cF \eta(\sqrt{\cD} {\boldsymbol{\nu}})
\, \ee^{ \frac{2 \pi i}{h} \langle {\mathbf x}, {\boldsymbol{\nu}}\rangle } 
= \sigma_{\boldsymbol{\alpha}}(\bx,\h,\cD) - \delta_{0|\boldsymbol{\alpha}|}
\, .
    \end{align*}
Hence \eqref{repre} leads to the representation of the quasi-interpolation error
   \begin{align*} %\label{fullexp}
(\cM_{h\sqrt{\cD}} &\,g) (\bx)-g(\bx) 
=(-h\sqrt{\cD})^N 
g_N(\bx)
+\sum _{ |\boldsymbol{\alpha}|  = 0} ^{N-1}
\frac{(h \sqrt{\cD})^{|\boldsymbol{\alpha}|}}{\boldsymbol{\alpha} !(2 \pi i)^{|{\boldsymbol{\alpha}}|}} 
\, \partial ^{\boldsymbol{\alpha}} g(\bx) \eps_{\boldsymbol{\alpha}} (\bx,\cD,\h) ,
    \end{align*}
which implies in particular the error estimate in $L^p(\R^n)$
 \begin{equation}  \begin{split} \label{estilp}
 \|g-&\cM_{h\sqrt{\cD}} \,g\|_{L^p}\! \le\! \\ &C_N (h\sqrt{\cD})^N|g|_{W_p^N}\!
+\!\sum _{ k = 0} ^{N-1} \frac{(h \sqrt{\cD})^k}{(2 \pi)^k}\! \sum _{ |\boldsymbol{\alpha}|  = k} 
\frac{\| \epsilon_{\boldsymbol{\alpha}}(\cdot,\cD,\h)\|_{L^\infty}
\| \partial ^{\boldsymbol{\alpha}} g\|_{L^p}}{\boldsymbol{\alpha}! } \, .
    \end{split}\end{equation}
Thus the quasi-interpolation error consists
of a term ensuring $\cO(h^N)$-convergence and
of the so-called saturation error, which, in general, does not converge to zero as $h \to 0$.
However,  due to the fast decay of $\partial^{\boldsymbol{\alpha}}\cF \eta$,  one can choose
$\cD$ large enough to ensure that 
   \begin{align*} 
\| \epsilon_{\boldsymbol{\alpha}}(\cdot,\cD,\h)\|_{L^\infty} \le 
\sum_{\boldsymbol{\nu} \in {\Z}^n\setminus \boldsymbol{0}}
|\partial^{\boldsymbol{\alpha}}\cF \eta(\sqrt{\cD} {\boldsymbol{\nu}})|<\e
    \end{align*}
for given small $\e > 0$.
In the examples below the saturation error is of the order $\cO(\ee^{-\pi^2 \cD})$,
which in the cases $\cD=2$ and $\cD=4$ 
is comparable to the single and double precision
arithmetics of modern computers. Therefore,
in numerical computations the saturation error can be neglected for appropriately chosen $\cD$.

\subsection{Approximation error for  $\cS g$}

From \eqref{appsol0} we see that
\begin{align} \label{semiS}
\cS_h g(\bx,t)= \cD^{-n/2} \sum_{\bm \in \Z^n} g(h \bm) \, 
\cS \h\Big(\frac{\bx-h\bm}{h\sqrt{\cD}},\frac{t}{h^2 \cD}\Big)  .
\end{align}
Hence, if $\cS \h(\bx,t)$ is known analytically, then \eqref{semiS}
is a very simple semi-analytic cubature of $\cS g$.
Of course, the cubature formula is computable only for a finite number
of  nonvanishing terms in \eqref{semiS}.
Therefore we assume that $g(\bx)$ and  $f(\bx,t)$  are compactly supported.  

The mapping properties \eqref{estip} and \eqref{striest} of $\cS$ 
and the quasi-interpolation error \eqref{estilp}
lead to estimates of the approximation error for $\cS g$.
In the following theorem
we use the notation
\begin{align*}
\|\nabla_k g\|_{L^{p}} = 
 \sum _{ |\boldsymbol{\alpha}|  = k} 
 \frac{\| \partial ^{\boldsymbol{\alpha}} g\|_{L^{p}}}{\boldsymbol{\alpha}! } \, .
\end{align*}

\begin{thm}\label{thm1}
Let $g \in W_{p}^N(\R^n)$, $1 \le p \le 2$, $N> n/p$, be  the initial value for the homogeneous 
Schr\"odinger equation \eqref{main0}.
For any $\e>0$ there exists $\cD > 0$
such that for $t>0$ the cubature formula \eqref{semiS}
approximates  the solution $u(\bx,t)$ in $L^{p'}(\R^n)$, $p'=p/(p-1)$, with
\begin{align*}
\|u(\cdot,t)- \cS_h g(\cdot,t)\|_{L^{p'}} \le 
\frac{C}{t^{n (1/2-1/p')}} \Big((h\sqrt{\cD})^N|g|_{W_{p}^N} 
+\e \sum _{ k = 0} ^{N-1} \frac{(h \sqrt{\cD})^k}{(2 \pi)^k} 
\|\nabla_k g\|_{L^{p}}\Big).
\end{align*}
Moreover, if $g \in W_{2}^N(\R^n)$, then the approximation on $\R^n \times \R$ 
with
$u_h(\bx,t)=\cS_h g(\bx,t)$ 
can be estimated in the mixed Lebesgue spaces $L^{r,q}(\R_+)$
for any  Schr\"odinger-admissible pairs $(q,r)$,
cf. \eqref{defadmissible},
by
\begin{align*}
\|u-u_h\|_{L^{r,q}(\R_+)} \le 
C \Big((h\sqrt{\cD})^N|g|_{W_{2}^N} 
+\e \sum _{ k = 0} ^{N-1} \frac{(h \sqrt{\cD})^k}{(2 \pi)^k} 
\|\nabla_k g\|_{L^{2}}\Big).
\end{align*}
\end{thm}

\begin{remark}
The analysis of the application of the integral operator
$\cS$ on the saturation error 
\begin{align*}
R_N(\bx)=\sum _{ |\boldsymbol{\alpha}|  = 0} ^{N-1}
\frac{( h \sqrt{\cD})^{|\boldsymbol{\alpha}|}}{\boldsymbol{\alpha} !(2 \pi i)^{|{\boldsymbol{\alpha}}|}} 
\, \partial ^{\boldsymbol{\alpha}} g(\bx) \, \eps_{\boldsymbol{\alpha}}(\bx,\h,\cD) 
\end{align*}
shows that  $|\cS R_N(\bx,t)|\to 0$ for fixed $(\bx,t)$
as $h \to 0$.
For example, if $(1+|\bx|^2)^{(N-1)/2}\\ \times  \!g(\bx)\in W_{1}^{N-1}(\R^n)$,
then $|\cS R_N(\bx,t)| = C t ^{-n/2} (\sqrt{\cD} h)^N \| (1+|\cdot|^2)^{(N-1)/2} g\|_{W_{1}^{N-1}}$ 
for all $|\bx| < 2 \pi t/h$.
\end{remark}

\subsection{Approximation error for  $\Pi f$}
We construct an approximation for $\Pi f$ in \eqref{sol1}
using the approximate quasi-interpolant 
\begin{equation}\label{fht}
\cN_{h\sqrt{\cD},\t\sqrt{\cD_0}} \, f(\bx,t) =\frac{1 }{\sqrt{\cD_0{\cD}^{n}}}
\sum_{\substack{\ell\in\Z\\\bm\in\Z^n}} 
f(h \bm,\t \ell)\, \psi
\Big(\frac{t-\t \ell}{\t \sqrt{\cD_0}}\Big) {\h} \Big(\frac{\bx-h \bm}{h \sqrt{\cD}}\Big) .
\end{equation}
Here $\t$, $h$, are the step sizes,
$\cD_0$ and $\cD$ are positive fixed
parameters; $\psi\in S(\R)$ and $\h\in S(\R^n)$ are the generating
functions, which belong to the Schwartz space $S$ of smooth and
rapidly decaying functions. If the generating functions $\psi$  and
$\h$
fulfills the moment condition \eqref{momcon} of order $N$, then
$\cN_{h\sqrt{\cD},\t\sqrt{\cD_0}} f$ approximates $f$ with the order
$\cO((h\sqrt{\cD}+\t \sqrt{\cD_0})^N)$ up to the saturation error in $L^p(\R^{n+1})$ if
$N>(n+1)/p$.
This is also true for mixed Lebesgue spaces, which are used
for the mapping properties  \eqref{striest} of $\Pi$.

More precisely, using the previously mentioned approach one can expand
\begin{equation*}\label{baserrNN}
(\cN_{h\sqrt{\cD},\t\sqrt{\cD_0}}  f)(\bx,t) = f_N(\bx,t)+ \cR_N f(\bx,t) \, ,
\end{equation*}
where
\begin{align*} 
&f_N(\bx,t) = \frac{(-1)^N N}{\sqrt{\cD_0{\cD}^{n}}} 
\sum_{k=0}^N \frac{(\t \sqrt{\cD_0})^k(h\sqrt{\cD})^{N-k}}{k! }\\& 
\times \sum _{ |\boldsymbol{\alpha}|  = N-k}
\frac{1}{\ba!}
\sum_{\substack{\ell\in\Z\\\bm\in\Z^n}} 
\Big(\frac{t-\t \ell}{\t \sqrt{\cD_0}}\Big)^k
\psi \Big(\frac{t-\t \ell}{\t \sqrt{\cD_0}}\Big)
\Big(\frac{\bx\!-\!h\bm}{h\sqrt{\cD}} \Big)^{\boldsymbol{\alpha}}\h \Big(
\frac{\bx\!-\!h\bm}{h\sqrt{\cD}} \Big) U_{k,\boldsymbol{\alpha}}(\bx,h \bm,t , \t \ell)
\end{align*}
with the notation
\begin{align*} 
U_{k,\boldsymbol{\alpha}}(\bx,\by,t , z)=
\itg _0 ^1 s^{N-1}  \partial_\bx ^{\boldsymbol{\alpha}} \partial_t^k f
     \big(s \bx + (1-s) \by,s t + (1-s)z \big)  \, ds \, ,
\end{align*}
and the function
\begin{align*} 
\cR_N f(\bx,t)=\sum_{k=0}^{N-1} \frac{(\t \sqrt{\cD_0})^k\, \sigma_{k}(t,\psi,\cD_0)}{k! \,(2 \pi i)^k}
\sum _{ |\boldsymbol{\alpha}|  = 0} ^{N-1-k} 
\frac{(h \sqrt{\cD})^{|\boldsymbol{\alpha}|}\, \sigma_{\boldsymbol{\alpha}}(\bx,\h,\cD)}
{\boldsymbol{\alpha}! \,(2 \pi i)^{|\boldsymbol{\alpha}|}}\, 
\partial_\bx ^{\boldsymbol{\alpha}} \partial_t^k f(\bx,t)\, .
\end{align*}
The moment conditions for $\psi$ and $\h$ obviously imply, that for $q,r \in [1,\infty]$ and  $f$ such that
the partial derivatives 
$ \partial_\bx ^{\boldsymbol{\alpha}} \partial_t^k  f\in L^{r,q}(\R)$
for all indices $|\boldsymbol{\alpha}|+k\le N$,
\begin{align*}
\|\cR_N f 
-f\|_{L^{r,q}} \le \|\sigma_{\boldsymbol{0}}(\cdot,\cD,\h)\|_{L^\infty}\sum_{k=0}^{N-1} \frac{(\t \sqrt{\cD_0})^k}{k! \,(2 \pi)^k} 
\, \|\epsilon_{k}(\cdot,\cD_0,\psi) \|_{L^\infty}
\|  \partial_t^k f\|_{L^{r,q}} \nonumber \\
+\|\sigma_{0}(\cdot,\cD_0,\psi) \|_{L^\infty} \sum_{k=0}^{N-1} \frac{(h \sqrt{\cD})^k}{(2 \pi)^k} 
\sum _{ |\boldsymbol{\alpha}|  = k}  \frac{ \|\eps_{\boldsymbol{\alpha}}(\cdot,\cD,\h) \|_{L^\infty}
}{ \boldsymbol{\alpha} !} \, \| \partial_\bx^{\boldsymbol{\alpha}} f\|_{L^{r,q}} \\
+\sum_{k=1}^{N-1} \frac{(\t \sqrt{\cD_0})^k\|\epsilon_{k}(\cdot,\cD_0,\psi) \|_{L^\infty}}{k! \,(2 \pi)^k}
\sum _{ |\boldsymbol{\alpha}|  = 1} ^{N-1-k} \frac{(h \sqrt{\cD})^{|\boldsymbol{\alpha}|}
\| \epsilon_{\boldsymbol{\alpha}}(\cdot,\cD,\h)\|_{L^\infty}}{\boldsymbol{\alpha} !(2 \pi)^{|\boldsymbol{\alpha}|}}
\|\partial_\bx ^{\boldsymbol{\alpha}} \partial_t^k f\|_{L^{r,q}} .  
\end{align*}
Moreover, if $N > (n+1)/\min(q,r)$,
then one can show similar to \cite[Lemma 2.29]{MSbook} that
\begin{align*} 
\|f_N\|_{L^{r,q}} \le C \sum_{k=0}^N \sum _{ |\boldsymbol{\alpha}|  = N-k} (\t \sqrt{\cD_0})^k 
(h\sqrt{\cD})^{N-k} \| \partial_t ^k \partial_\bx ^{\boldsymbol{\alpha}}f\|_{L^{r,q}} \, .
\end{align*}

Now we are in the position to study the approximation  of $\Pi f$
defined by 
\begin{align*} 
\Pi_{h,\t} f(\bx,t) \!=\!\Pi(\cN_{h\sqrt{\cD},\t\sqrt{\cD_0}}  f)(\bx,t) =\!
-i\! \itg_0^{t} \!\!ds\! \itg_{\R^n} \!\cK(\bx-\by,t-s) \,\cN_{h\sqrt{\cD},\t\sqrt{\cD_0}}  f (\by,s) \,  d\by .
\end{align*}
The difference between $\Pi f$ and $\Pi_{h,\t} f$ 
can be estimated by \eqref{striest} and the quasi-interpolation error 
$\| f- \cN_{h\sqrt{\cD},\t\sqrt{\cD_0}}  f\|_{L^{r',q'}}$.
However, the corresponding error estimate is proved for sufficiently smooth functions on $\R^{n+1}$,
whereas the right-hand side $f$ of \eqref{inhom} is given only on $\R^n\times\R_+$.
Therefore we extend $f(\cdot,t): \R_+ \to W_r^N(\R^n)$ with preserved
smoothness to a function $f(\cdot,t): \R \to W_r^N(\R^n)$.
This can be done, for example, by using Hestenes
reflection principle (cf. \cite{He} and \eqref{w}).  The  extended function, again denoted by
$f$, is compactly supported and  
retains the smoothness of $f|_{\R^n\times\R_+}$. Then we get
\begin{align} \nonumber
\Pi_{h,\t}& f(\bx,t)\\
&=\frac{-i}{\sqrt{\cD_0{\cD}^{n}}}
\sum_{\substack{\ell\in\Z\\\bm\in\Z^n}} f(h \bm,\t \ell) 
\itg_0^{t}\psi \Big(\frac{t-\t \ell-s}{\t \sqrt{\cD_0}}\Big)\itg_{\R^n} 
\frac{\ee^{ \, i \, |\bx-\by|^2/(4s)}}{(4 \pi i s)^{n/2}}
 {\h} \Big(\frac{\by-h \bm}{h \sqrt{\cD}}\Big) d\by ds\nonumber \\
&= \frac{-i }{\sqrt{\cD_0 \,\cD^{n}}}\sum_{\substack{\ell\in\Z\\\bm\in\Z^n}} f(h \bm,\t \ell) 
\itg_0^{t}\psi \Big(\frac{t-\t \ell -s}{\t \sqrt{\cD_0}}\Big)
\cS \h \Big(\frac{\bx-h \bm}{h \sqrt{\cD}},\frac{s}{h^2{\cD}}\Big)\,ds \,  \label{semiPi}
\end{align}
where  $f(h \bm,\t \ell)$ for $\ell < 0$ 
in \eqref{semiPi} are understood as values of the extended function.

\begin{thm}\label{thmPi}
Let $(q,r)$ be a Schr\"odinger-admissible pair and $N>(n+1)/\min(q',r')$, $q'=q/(q-1)$, $r'=r/(r-1)$.
Suppose that the right-hand side $f$ of the
inhomogeneous Schr\"odinger equation \eqref{inhom} satisfies 
$\partial_t ^k \partial_\bx ^{\boldsymbol{\alpha}}f \in  L^{r',q'}(\R_+)$
for all $0\le k+ |\boldsymbol{\alpha}| \le N$.
Then there exist a constant $C$ and
for any $\e>0$  parameters  $\cD_0,\cD > 0$, not depending on $f$, 
such that the cubature formula \eqref{semiPi}
provides the approximation estimate
\begin{align*} 
\| \Pi f-\Pi_{h,\t}   f\|_{L^{r,q}(\R_+)}&\le  \, C 
\sum_{k=0}^N \sum _{ |\boldsymbol{\alpha}|  = N-k} (\t \sqrt{\cD_0})^k 
(h\sqrt{\cD})^{N-k} \| \partial_t ^k \partial_\bx ^{\boldsymbol{\alpha}}f\|_{L^{r',q'}(\R_+)}\\
& +\e \sum _{ k= 0} ^{N-1} \sum _{ j = 0} ^{N-1-k} 
\frac{(\t \sqrt{\cD_0})^k(h \sqrt{\cD})^j}{(2 \pi)^{k+j}} 
\sum _{ |\boldsymbol{\alpha}|  = j}\| \partial_t ^k \partial_\bx ^{\boldsymbol{\alpha}}f\|_{L^{r',q'}(\R_+)}
\, .
\end{align*}
\end{thm}

\section{Cubature formulas}\label{sec3}
\setcounter{equation}{0}

\subsection{Approximation of $\cS g$ }

For different basis functions $\h$ the integrals on the right in \eqref{appsol0} allow analytic
representations.
For example, let for $N=2M$
\begin{equation*} \label{eta2M}
          \eta(\bx)=\eta_{N}(\mathbf{x}) =\pi^{-n/2} \, L_{M-1}^{(n/2)}(|\mathbf{x}|^2) 
        \, \ee^{-|\mathbf{x}|^2} 
\end{equation*}
with the generalized Laguerre polynomials
   \begin{equation*}
L_{k}^{(\gamma)}(y)=\frac{\ee^{\,y} y^{-\gamma}}{k!} \, \Big(
{\frac{d}{dy}}\Big)^{k} \!
\left(\ee^{\,-y} y^{k+\gamma}\right), \quad \gamma > -{\rm 1} \, .
    \end{equation*}
By using the relation (\cite[Theorem 3.5]{MSbook})
\[
\eta_{N} (\bx)=
\pi^{-n/2} \sum_{j=0}^{M-1}
        \frac{(-1)^j}   {j!\, 4^j }
\Delta^j  \ee^{\, - |\mathbf{x}|^2} 
\]
we obtain the formula
\begin{align*}
\frac{1}{(\pi i t)^{n/2}}  \itg_{\R^n}\ee^{i |\bx-\by|^2/t}\h_N (\by)d\by 
&=\frac{1}{\pi^n (it)^{n/2}}
\sum_{j=0}^{M-1}\frac{(-1)^j}   {j! \, 4^j}
\Delta^j \itg_{\R^n} \ee^{\, i \, |\mathbf{x}-\mathbf{y}|^2/t} 
 \ee^{- |\mathbf{y}|^2 } \, d \mathbf{y} \\
&= \frac{1}{\pi^{n/2} (1 +it)^{n/2}}
\sum_{j=0}^{M-1}\frac{(-1)^j}   {j! \, 4^j }
\Delta^j \ee^{\, - |\mathbf{x}|^2/(1 +it)} \, .
\end{align*}
In view of (\cite[(3.15)]{MSbook})
\[
\Delta^j \ee^{\, - |\mathbf{x}|^2}=(-1)^j \, j! \, 4^j \ee^{\, - |\mathbf{x}|^2}L_j^{(n/2-1)}( |\bx|^2) 
\]
an approximate solution of \eqref{main0} is given by the analytic formula
\begin{align*}%\label{Nsol}
&\cS_h g (\mathbf{x},t)= \cD^{-n/2} \sum_{\bm \in \Z^n} g (h{\mathbf m})
\Phi_{N}\Big(\frac{\bx-h \bm}{h \sqrt{\cD}},\frac{4\,t}{h^2\cD}\Big)
\end{align*}
with
\[
\Phi_{N}(\bx,t)=\frac{\ee^{-|\bx|^2/(1+i\, t)}}{\pi^{n/2}(1+i\, t)^{n/2}} 
\sum_{j=0}^{M-1}  \frac{1}{(1+i\, t)^j} \, L_j^{(n/2-1)}\Big(\frac{|\bx|^2}{1+i\, t}\Big).
\]

Another approximation of order $N=2M$ of the initial function $g$
can be derived by the quasi-interpolant
\[
g_h(\bx)=\cD^{-n/2} \sum_{\bm \in \Z^n} g(h \bm) \widetilde{\h_N}\left(\frac{\bx-h\bm}{h\sqrt{\cD}}\right)
\]
with a basis function in tensor product form
\begin{equation}
\label{tensorbasis}
\widetilde{\h_N}(\bx)=\prod_{j=1}^n \chi_{2M}(x_j);
\quad
{\chi}_{2M}(x_j)=\frac{(-1)^{M-1}}{2^{2M-1}\sqrt{\pi} (M-1)!}\frac{H_{2M-1}(x_j) \ee^{-x_j^2}}{x_j}.
\end{equation}
$H_k$ are the Hermite polynomials
\begin{equation*}
H_k(x)=(-1)^k \ee^{x^2} \Big( \frac{d}{dx}\Big)^k \ee^{-x^2}.
\end{equation*}
Then $u(\bx,t)$ in \eqref{sol0} is approximated by 
\begin{equation} \label{tensorsol}
 \cS_h  g(\bx,t)=\cD^{-n/2}\sum_{\bm\in\Z^n}g(h \bm) \widetilde{\Phi_{N}}\Big(\frac{\bx-h \bm}{h \sqrt{\cD}},\frac{4\,t}{h^2\cD}\Big)
\end{equation}
with
\begin{equation}\label{PhiNbis}
\widetilde{\Phi_{N}}(\bx,t)= \prod_{j=1}^n \phi_{2M}(x_j,t)=\prod_{j=1}^n
\frac{1}{(\pi\,i\,t )^{1/2}}
 \itg_{\R} \ee^{\,i (x_j-y)^2/t} \chi_{2M}(y) \, dy \, , \quad t\neq 0.
\end{equation}
From  the representation \cite[(3.9) and (3.6)]{MSbook}
\begin{equation*}
 \chi_{2M}(y)=\frac{1}{\sqrt{\pi}} \sum_{\ell=0}^{M-1} \frac{(-1)^\ell}{\ell! \,4^\ell}
\frac{\partial^{2\ell}}{\partial y^{2\ell}}\ee^{-{y^2}}
\end{equation*}
we obtain
\begin{align*}
\phi_{2M}(x_j,t)
&= \frac{1}{\sqrt{\pi}}
 \sum_{\ell=0}^{M-1} \frac{(-1)^\ell}{\ell! \, 4^\ell} \frac{\partial^{2\ell}}{\partial x_j^{2\ell}}
 \frac{\ee^{-x_j^2/(1+i\,t)}}{(1+i\, t)^{1/2}}\\& =
\frac{1}{\sqrt{\pi}}\sum_{\ell=0}^{M-1}  \frac{\ee^{-x_j^2/(1+i\,t)}}{(1+i\, t)^{\ell+1/2}} \, L_\ell^{(-1/2)}\Big(\frac{x_j^2}{1+i\, t}\Big),
\end{align*}
which in view of
\[
L_\ell^{(-1/2)}(y^2)= \frac{(-1)^\ell}{\ell! \, 4^\ell}H_{2\ell}(y)
\] 
gives
\begin{equation}\label{phiN}
\widetilde{\Phi_{N}}(\bx,t)=\frac{\ee^{-|\bx|^2/(1+i\, t)}}{\pi^{n/2}(1+i\, t)^{n/2}} \prod_{j=1}^n
\sum_{\ell=0}^{M-1} \frac{(-1)^\ell}{\ell! \, 4^\ell} \frac{1}{(1+i\, t)^\ell} H_{2\ell}\Big(\frac{x_j}{\sqrt{1+i\, t}}\Big).
\end{equation}

Note that the computation of the approximate solution with the summation \eqref{tensorsol}
is very efficient if the function $g(\bx)$ 
allows a separated
representation; that is, within a prescribed accuracy, it can be
represented as sum of products of univariate functions
\begin{equation}\label{sep}
g(\bx)=
\sum_{p=1}^P \alpha_p \prod_{j=1}^n g^{(p)}_j (x_j)+\cO(\e)\, .
\end{equation}
Then the products of one-dimensional sums 
\begin{align*}
 \cS_h g(\bx,t)\approx \frac{h^n}{\pi^{n/2}}&\sum_{p=1}^P \alpha_p \prod_{j=1}^n \sum_{m_j\in\Z}g^{(p)}_j (hm_j) 
\phi_{2M}\left(\frac{x_j-hm_j}{h \sqrt{\cD}},\frac{4\,t}{h^2\cD}\right) \\
=\frac{h^n}{\pi^{n/2}}&\sum_{p=1}^P \alpha_p \prod_{j=1}^n \sum_{m_j\in\Z}g^{(p)}_j (hm_j)\\
&\times 
\sum_{\ell=0}^{M-1} \frac{(-1)^\ell}{\ell! \,4^\ell} 
\frac{\ee^{-(x_j-hm_j)^2/(h^2\cD+4i\, t)}}{(h^2\cD+4i\, t)^{\ell+1/2}} H_{2\ell}\left(\frac{x_j-hm_j}{\sqrt{h^2\cD+4i\, t}}\right)
\end{align*}
provide the approximation of the $n$-dimensional
initial value problem for the  Schr\"odinger equation \eqref{main0}.

\subsection{Approximation of $\Pi f$ in  \eqref{sol1}} \label{subsec32}
Let us assume in \eqref{fht}
\[
\psi(t)=\chi_{2M}(t),\qquad \h(\bx)=\prod_{j=1}^n \chi_{2M}(x_j).
\]
Then $\Pi  f$  is approximated by
\begin{equation}\label{tensorsol2}
\begin{split}
&\Pi_{h,\t}  f(\bx,t)=\Pi (\cN_{h\sqrt{\cD},\t\sqrt{\cD_0}} f)(\bx,t)\\
&=-\frac{i}{\sqrt{\cD_0 \cD^{n}}} \sum_{\substack
{\ell\in\Z\\\bm\in\Z^n}} f(h \bm,\t \ell)\,   \itg_0^t 
\chi_{2M}\Big(\frac{s-\t \ell}{\t\sqrt{\cD_0}}\Big)
\, \widetilde{\Phi_{N}}\Big(\frac{\bx-h \bm}{h \sqrt{\cD}},\frac{4\,(t-s)}{h^2\cD}\Big) ds.
\end{split}
\end{equation}

The approximation of $\Pi f$ requires the computation of a certain
number of one-dimensional integrals where, for \eqref{phiN}, the
integrands allow separated representations.  Suppose that also
$f(\bx,t)$ allows a separated representation, that is
\begin{equation}\label{sep2}
{f}(\bx,t)=
\sum_{p=1}^P \beta_p \prod_{j=1}^n f^{(p)}_j (x_j, t)+\cO(\e),
\end{equation}
then, for \eqref{tensorsol2} and \eqref{PhiNbis},
\begin{align}\label{Sigmath}
\Pi_{h,\t}   f(\bx,t)\approx
\frac{-i}{\sqrt{\cD_0 \cD^{n}}} \sum_{\ell\in\Z} \sum_{p=1}^P \beta_p\itg_0^t \chi_{2M}\left(\frac{s-\t \ell}{\tau\sqrt{\cD_0}}\right)  \prod_{j=1}^n T_j^{(p)}(x_j,\frac{4\,(t-s)}{h^2\cD},\t \ell)ds
\end{align}
where
\begin{equation*}\label{TN}
T_j^{(p)}(x,t,\t \ell)=\sum_{m\in\Z } f^{(p)}_j (hm, \t \ell )\phi_N(\frac{x-hm}{h\sqrt{\cD}},t)\,.
\end{equation*}

An accurate quadrature rule of the one-dimensional integrals in
\eqref{Sigmath} provides a separated representation of $\Pi_{h,\t} f$
(this is described in detail in section \ref{sec4}). Then the
numerical computation of $\Pi_{h,\t} f$ does not require to perform
$n$-dimensional integrals and sums but only one-dimensional
operations, which leads to a considerable reduction of computing
resources, and gives the possibility to treat real world problems.

\section{Schr\"odinger equation over hyper-rectangles}\label{sec4}
\setcounter{equation}{0}

  Assume now that  $g$ in \eqref{initcond} is supported in a hyper-rectangle $[{\bf P,Q}]=\{\bx\in\R^n: P_j\leq x_j\leq Q_j, j=1,...,n\}$ and $g\in C^N([{\bf P,Q}])$. Hence 
 \begin{equation}\label{sol3}
 \cS  g(\bx,t)= \frac{1}{(4 \pi i t)^{n/2}}  \itg_{[{\bf P,Q}]} \ee^{i |\bx-\by|^2/(4t)} g(\by)d\by,\quad \bx\in\R^n,\quad t\in \R
\end{equation}
provides the solution of \eqref{main0} with the data \eqref{initcond}.

The direct application of the method described in Section \ref{sec3}  does not give good approximations because the sum
\[
\cD^{-n/2} \sum_{h \bm \in [\bP,\bQ]} g(h \bm) \h\left(\frac{\bx-h\bm}{h\sqrt{\cD}}\right)
\]
approximates $g$ only in a subdomain of $[\bP,\bQ]$. To overcome this
difficulty  we extend $g$  by using the Hestenes
reflection principle into a larger domain with preserved
smoothness.  
If $\widetilde g$ is the
extension of $g$, since $\h$ is of rapid decay, one can fix $r>0$ such
that the quasi-interpolant
\begin{equation*}\label{gh}
\cD^{-n/2} \sum_{h\bm \in \O_{rh}} \widetilde{g}(h \bm) \h\left(\frac{\bx-h\bm}{h\sqrt{\cD}}\right)
\end{equation*}
approximates $g$ in $[\bP,\bQ]$ with the same error estimate of \eqref{appN}.  Here  $\O_{rh}=\prod_{j=1}^n I_j$, $I_j=(P_j-rh\sqrt{\cD},Q_j+rh\sqrt{\cD})$.
 
In the following we consider the basis functions in tensor product form \eqref{tensorbasis}.
Then the sum
\[
 \cS_h ^{[\bP,\bQ]}g(\bx,t)=\cD^{-n/2}\sum_{\bm\in \O_{rh}}\widetilde{g}(h \bm)
\widetilde \Phi_{2M}^{[{\bf P_\bm,Q_\bm}]}(\frac{\bx-h \bm}{h \sqrt{\cD}},\frac{4\,t}{h^2\cD}),\quad t\neq 0
 \]
with $\bP_\bm=(\bP-h\bm)/(h \sqrt{\cD})$, $\bQ_\bm=(\bQ-h\bm)/(h \sqrt{\cD})$, and 
\begin{equation*}
\begin{split}
\widetilde\Phi_{2M}^{[{\bf P,Q}]}(\bx,t)=&\frac{1}{(\pi\,i\,t )^{n/2}} 
\itg_{[{\bf P,Q}]} \ee^{i\, |\by-\bx|^2/t} \widetilde{\eta_{2M}}(\by)d\by\\
=&\prod_{j=1}^n\frac{1}{(\pi\,i\,t )^{1/2}} \itg_{P_j}^{Q_j} 
\ee^{i\, (y_j-x_j)^2/t} \chi_{2M}(y_j)dy_j,
\end{split}
\end{equation*}
provides an approximation of  \eqref{sol3} with the error estimate obtained in 
Theorem \ref{thm1}. $\widetilde\Phi_{2M}^{[{\bf P,Q}]}(\bx,4\,t)$ gives the solution of the initial problem
\begin{equation}\label{init2}
i\partial_t v+\D_\bx v=0,\quad
v(\bx,0)=\prod_{j=1}^n I_{(P_j,Q_j)}(x_j)\chi_{2M}(x_j) , \quad \bx\in \R^n,\quad t\in \R.
\end{equation}
Here $I_{(P_j,Q_j)}$ is the characteristic function of the interval $(P_j,Q_j)$.
In \cite[Theorem 3.1]{LMS5} we prove
\begin{thm}\label{thm3}
The solution of the initial value problem \eqref{init2} in $\R^n$ can be expressed by the tensor product
\[
v(\bx,t)=\prod_{j=1}^n\big(\Psi_M(x_j,4 \,t , P_j)-\Psi_M(x_j,4 \,t , Q_j)\big),
\]
where
\begin{equation}\label{Psi}
\Psi_M(x,t,y)=\frac{1}{2\sqrt{\pi}}
\ee^{- x^2/(1+i\,t)}
 \bigg(\erfc(F(x,i\,t,y)) \cP_M(x,i\,t) 
-  \frac{\ee^{-F^2(x,i\,t,y)}}{\sqrt{\pi} } \cQ_M(x,i\,t,y)
\bigg)
\end{equation}
with the complementary error function $\erfc$, the argument function
\begin{equation}\label{defF}
F(x,t,y)=\sqrt{\frac{t+1}{t}} \Big(y-\frac{x}{t+1}\Big),
\end{equation}
and $\cP_M$, $\cQ_M$ are polynomials in $x$ of degree $2M-2$ and $2M-3$, respectively:
\begin{equation*}\label{QM}
\begin{split}
\cP_M(x,t)&=
\sum_{s=0}^{M-1} \frac{(-1)^s}{s! 4^s } \frac{1}{(1+  t)^{s+1/2}} H_{2s}\left(\frac{ x}{\sqrt{1+  t}}\right)\,;
\\
\cQ_1(x,t,y)&=0, \\%\nonumber \\
\cQ_M(x,t,y) &=2
\sum_{k=1}^{M-1} \frac{(-1)^k}   {k! \, 4^k}
  \sum_{\ell=1}^{2k}
\frac{(-1)^{\ell}}{t^{\ell/2}}
\bigg( H_{2k-\ell} (y)
  H_{\ell-1} \Big( \frac{y-x}{\sqrt{t}}\Big) \\
&\hskip60pt -
\Big(\hskip-1pt\begin{array}{c}2k\\\ell\end{array}\hskip-1pt\Big)
H_{2k-\ell}\Big(\frac{x}{\sqrt{1+t}}\Big)
\frac{H_{\ell-1}\big(F(t, x, y)\big)}{(1+t)^{k+1/2}}
\bigg)\,, \> M>1.
\end{split}
\end{equation*}

\end{thm}

From Theorem \ref{thm3} we deduce the following semi-analytic cubature formula  for \eqref{sol3} with the error $\cO((h\sqrt{\cD})^{2M})$
\begin{equation*}
\begin{split}
 \cS  g(\bx,t)&\approx \frac{1}{\cD^{n/2}}\sum_{h\bm\in\Omega_{rh}}\widetilde{g}(h \bm) \\
\times & \prod_{j=1}^n\left(
 \Psi_M\big(\frac{x_j-hm_j}{h\sqrt{\cD}},\frac{4\,t}{h^2\cD},\frac{P_j-hm_j}{h\sqrt{\cD}}\big)
-\Psi_M\big(\frac{x_j-hm_j}{h\sqrt{\cD}},\frac{4\,t}{h^2\cD},\frac{Q_j-hm_j}{h\sqrt{\cD}}\big)
 \right).
 \end{split}
  \end{equation*}
  If $\widetilde g$ allows a separated representation \eqref{sep} we
derive that, at the points of the uniform grid $\{h\bk,\t s \}$, the
$n-$ dimensional integral \eqref{sol3} is approximated by the product
of one-dimensional sums
\[
\frac{1}{\cD^{n/2}}\sum_{p=1}^P \a_p \prod_{j=1}^n S_j^{(p)}(h\bk,\t s)
\]
where
\begin{equation*}
\begin{split}
&S_j^{(p)}(h\bk,t)=\\ &\sum_{hm\in I_j }{g_j^{(p)}}(h m) 
 \Big(
 \Psi_M\big(\frac{k_j-m}{\sqrt{\cD}},\frac{4\,t}{h^2\cD},\frac{P_j-hm}{h\sqrt{\cD}}\big)
-\Psi_M\big(\frac{k_j-m}{\sqrt{\cD}},\frac{4\,t}{h^2\cD},\frac{Q_j-hm}{h\sqrt{\cD}}\big)
 \Big).
\end{split}
  \end{equation*}

Suppose now that the source term $f(\bx,t)$ in \eqref{inhom} is
supported with respect to $\bx$ in the hyper-rectangle $[{\bf
P,Q}]=\{\bx\in\R^n: P_j\leq x_j\leq Q_j, j=1,...,n\}$ and $f\in
C^N([{\bf P,Q}]\times \R)$. Then, from \eqref{sol1},
\begin{equation}\label{sol2}
\Pi  f(\bx,t)=-i \itg_0^{t} \frac{ds}{(4 \pi i s)^{n/2}} \itg_{[\bP,\bQ]} \ee^{i|\bx-\by|^2/(4\,s)} f(\by,t-s)d\by
\end{equation}
provides the solution of \eqref{inhom} with null initial data. We
extend $f(\cdot,t)$ outside $[\bP,\bQ]$ with preserved smoothness and
denote by $\widetilde f$ its extension. Due to the rapid decay of
the generating function $\widetilde \h$, one can fix $r$ and $r_0$,
positive parameters, such that the quasi-interpolant
\begin{equation*} 
\cN_{h\sqrt{\cD},\t\sqrt{\cD_0}}^{(r,r_0)} f(\bx,t) =
\frac{1}{\cD_0^{1/2} \cD^{n/2}}\!\!\sum_{\substack
{{h\bm \in \O_{rh}}\\ {\t\ell \in \widetilde\O_{r_0\t}    } }   }\! \!\!\widetilde{f}(h
\bm,\t \ell)\, 
\chi_{2M}\left(\frac{t-\t \ell}{\t \sqrt{\cD_0}}\right)\prod_{j=1}^n \chi_{2M} \!
\left(\frac{x_j-h m_j}{h \sqrt{\cD}}\right)
\end{equation*}
approximates $f$ for all $\bx\in[\bP,\bQ]$ and for all $t\in[-T,T]$, $T>0$, 
with order $\cO((h\sqrt{\cD}+\t \sqrt{\cD_0})^N)$. Here 
$ \widetilde\O_{r_0\t}=(-T-r_0 \t \sqrt{\cD},T+r_0 \t \sqrt{\cD})$. Hence
\begin{align*}%\label{Sigmath2}
\Pi_{h,\t}   f(\bx,t)=
\frac{-i}{\cD_0^{1/2} \cD^{n/2}}\sum_{\substack
{{h\bm \in \O_{rh}}\\ {\t\ell \in \widetilde\O_{r_0\t}    } }   } \widetilde{f}(h
\bm,\t \ell) K_{2M}(\bx,t, h \bm, \t \ell) ,
\end{align*}
where
\begin{equation}\label{K2M}
 K_{2M}(\bx,t, h \bm, \t \ell) = \itg_0^t \chi_{2M}\Big(\frac{t-s-\t \ell}{\t \sqrt{\cD_0}}\Big) 
\widetilde \Phi_{2M}^{[{\bf P_\bm,Q_\bm}]}\Big(\frac{\bx-h \bm}{h \sqrt{\cD}},\frac{4\,s}{h^2\cD}\Big) ds
\end{equation}
and
\[
\widetilde\Phi_{2M}^{[{\bf P,Q}]}(\bx,t)=\prod_{j=1}^n\big(\Psi_M(x_j,t , P_j)-\Psi_M(x_j,t , Q_j)\big).
\]
The integrals in \eqref{K2M} cannot be taken analytically. Therefore
we use an efficient quadrature based on the classical trapezoidal
rule, which is exponentially converging for rapidly decaying smooth
functions on the real line. Making the substitution introduced in
\cite{Mo}
\begin{equation*}
s=t \varphi(\xi),\qquad  \varphi(\xi)=\frac{1}{2} \Big(1+\tanh \big( \frac{a \pi}{2} \sinh \x \big) \Big)
= \frac{1}{1+\ee^{-a \pi \sinh \x}} \, ,
\end{equation*}
with certain positive constant $a$, $K_{2M}$ transforms to
\[
 K_{2M}(\bx,t, h \bm, \t \ell)\! =\!\frac{\pi a t}{2}\!\!\! \itg_{-\infty}^\infty \!\!\!
  \chi_{2M}\Big(\frac{t(1- \varphi(\xi))-\t \ell}{\t \sqrt{\cD_0}}\Big) 
\widetilde \Phi_{2M}^{[{\bf P_\bm,Q_\bm}]}
\Big(\frac{\bx-h \bm}{h \sqrt{\cD}},\frac{4\,t \varphi(\xi)}{h^2\cD}\Big) \o(\xi)d\xi
\]
where we denote
\[
\o(\xi)=\frac{\cosh \xi}{1+\cosh(a \pi \sinh \xi)}\,.
\]
The trapezoidal rule with step size $\k$ gives for sufficiently large $R\in\N$
  \begin{equation*}
  \begin{split}
K_{2M}(\bx,t, h \bm, \t \ell)& \approx \\ \frac{\pi a t \k}{2} \sum_{q=-R}^R &
  \chi_{2M}\Big(\frac{t(1-\varphi(\k q))-\t \ell}{\t \sqrt{\cD_0}}\Big) 
\widetilde \Phi_{2M}^{[{\bf P_\bm,Q_\bm}]}
\Big(\frac{\bx-h \bm}{h \sqrt{\cD}},\frac{4\,t \varphi(\k q)}{h^2\cD}\Big) \o(\k q).
\end{split}
  \end{equation*}
In general the number $R$ in the quadrature rule depends on $t$
because the integrand depends on $t$.  However, 
since $|\varphi(\xi)|\leq 1$ and $| \chi_{2M}(\xi)|\leq m_1$, $\forall \xi\in \R$,
it follows, that for $t\in[0,T]$, $\forall \bx\in[\bP,\bQ]$ and
$h\bm\in \Omega_{rh}$ 
\[
|\widetilde \Phi_{2M}^{[{\bf P_\bm,Q_\bm}]}
\Big(\frac{\bx-h \bm}{h \sqrt{\cD}},\frac{4\,t \varphi(\xi)}{h^2\cD}\Big)|\leq m_2\,.
\]

Hence $R$ can be choosen (independent of $t$) such that the
trapezoidal rule with step $\kappa$ gives an accurate approximation of
the integral over $\R$ of $\omega(\xi)$, which decays doubly exponentially as $|\xi|\to
\infty$.

We obtain that, at the points of the uniform grid $\{h\bk, \t s\}$,
the $n-$dimensional integral \eqref{sol2} is approximated by
\begin{align*}
\Pi_{h,\t}   f(h \bk,\t s)& \approx 
\frac{-i\,\pi a \t s \k}{2\cD_0^{1/2} \cD^{n/2}} 
\sum_{\substack{h\bm \in \O_{rh} \\ \t\ell \in \widetilde\O_{r_0\t} }} \widetilde{f}(h\bm,\t \ell)\\
& \times \sum_{q=-R}^R
  \chi_{2M}\Big(\frac{s(1-\varphi(\k q))-\ell}{\sqrt{\cD_0}}\Big) 
\widetilde \Phi_{2M}^{[{\bf P_\bm,Q_\bm}]}
\Big(\frac{\bk-\bm}{\sqrt{\cD}},\frac{4\,\t \ell\varphi(\k q)}{h^2\cD}\Big) \o(\k q).
\end{align*}
If $\widetilde f$ allows a separated representation \eqref{sep2} we
get the efficient high-order approximation
\begin{equation}\label{Sigmath3}
\begin{split}
\Pi_{h,\t}   f(h \bk,\t s)&\approx
\frac{-i\,\pi a \t s \k}{2\cD_0^{1/2} \cD^{n/2}}  \sum_{q=-R}^R
\o(\k q) \\ &\times\sum_{\t\ell \in \widetilde\O_{r_0\t} } \chi_{2M}
\Big(\frac{s(1-\varphi(\k q))-\ell}{\sqrt{\cD_0}}\Big) 
 \sum_{p=1}^P \beta_p \prod_{j=1}^n T^{(p)}_j(k_j,\t s,\t \ell,\k q)\, ,
\end{split}
\end{equation}
 where
\begin{equation}\label{Sigmath4}
\begin{split}
T^{(p)}_j(k_j,\t s,\t \ell,\k q)&= \sum_{
{{hm_j \in I_j} }   } f^{(p)}_j (h m_j, \t \ell) 
\left(\Psi_M\Big(\frac{k_j-m_j}{\sqrt{\cD}},
\frac{4\t s\varphi(\k q)}{h^2\cD},\frac{P_j-hm_j}{\sqrt{\cD}}\Big) \right. \\
& \left. -  \Psi_M\Big(\frac{k_j-m_j}{\sqrt{\cD}},
\frac{4\,\t s \varphi(\k q)}{h^2\cD},\frac{Q_j-hm_j}{\sqrt{\cD}}\Big) \right).
\end{split}
\end{equation}
We obtain that, if $f$ has the form \eqref{sep2}, then the approximation of
the potential \eqref{sol2} requires us to compute $2 R P n$
one-dimensional sums. Thus, if $n>1$, the computational time scales
linearly in the space dimension $n$.

For an efficient implementation of $\Psi_M$ we  express  
$\erfc$ in \eqref{Psi}    
with the Faddeeva or scaled 
complementary error function 
$W(z)=\ee^{-z^2} \erfc(-i z)$
(cf. \cite[7.1.3]{AS})
and write
\begin{align*}
\Psi_M(x,t,y)&=\frac{\ee^{-x^2/(1+it)-F^2(x,it,y)}}{2\sqrt{\pi}}
\left(  W\big(i F(x,it,y)\big)\cP_M(x,it)-\frac{\cQ_M(x,it,y)}{\sqrt{\pi} } 
\right) \\
&=\frac{\ee^{-y^2+i(y-x)^2/t}}{2\sqrt{\pi}}
\left( W\big(i F(x,it,y)\big)\cP_M(x,it) -\frac{\cQ_M(x,it,y)}{\sqrt{\pi} } 
\right),
\end{align*}
where $F(x,it,y)$ is defined by \eqref{defF}. 
Efficient implementations of  double precision computations of $W(z)$ 
 are available
if the imaginary part of the argument is nonnegative.
Otherwise, for $\im z < 0$,
overflow problems can occur,
which can be seen from the relation $W(z)=2\ee^{-z^2}-W(-z)$ (cf. \cite[7.1.11]{AS}).
But this helps to derive a stable formula also for $\im(i F(x,it,y))=\re F(x,it,y) < 0$,
since
\begin{align*}
\frac{\ee^{-x^2/(1+it)-F^2(x,it,y)}}{2} &W(iF(x,it,y))\\
 &=\frac{\ee^{-x^2/(1+it)-F^2(x,it,y)}}{2}\Big(2 \ee^{F^2(x,it,y)}-W\big(-i F(x,it,y)\big)\Big)\\
&=\ee^{-x^2/(1+it)}-\frac{\ee^{-y^2+i(y-x)^2/t}}{2}W\big(-i F(x,it,y)\big) \, .
\end{align*}
Thus we get the efficient formula
\begin{equation}\label{PsiNew}
\begin{split}
&\Psi_M(x,t,y)= \, -\frac{\ee^{-y^2+i(y-x)^2/t}}{2\sqrt{\pi}}
\frac{\cQ_M(x,it,y)}{\sqrt{\pi}} +\\
&\left\{
\begin{aligned}
& \ee^{-y^2+i(y-x)^2/t} \, 
W\big(i F(x,it,y)\big)\frac{\cP_M(x,it)}{2\sqrt{\pi}}  , &\re F(x,it,y) \geq 0 ,\\
&\Big(2\,{\ee^{-x^2/(1+it)}}-{\ee^{-y^2+i(y-x)^2/t}} \, W\big(-i F(x,it,y)\big)\Big)\frac{\cP_M(x,it)}{2\sqrt{\pi}},
& \re F(x,it,y)<0 .
\end{aligned}
\right.
\end{split}
\end{equation}

\section{Numerical Tests}\label{sec5}
 \setcounter{equation}{0}
 
 In this section we present some numerical results. First we verify
numerically the accuracy and the convergence order of the proposed
method for the inhomogeneous Schr\"odinger equation \eqref{inhom}
with null initial data and then for the initial value problem
\eqref{main0}-\eqref{initcond}.  Finally, in Figures
\ref{figre30}-\ref{figresin}, we depict the evolution of $u(\bx,t)$
under the two-dimensional equation
\eqref{main0}  for different initial values.
\subsection{Inhomogeneous Schr\"odinger equation}

We consider the Cauchy problem
\begin{equation}\label{main00}
i \frac{\partial u}{\partial t} +\D_\bx u=f(\bx,t),\quad u(\bx,0)=0\quad \bx\in\R^n
\end{equation}
for right-hand sides
\begin{equation}\label{fxt}
f(\bx,t)=\left( i \frac{\partial }{\partial t} +\D_\bx \right) \prod_{j=1}^n w(x_j) v(t)
\end{equation}
with  $\supp w\subset [-1,1]$. If $w(\pm 1)=w'(\pm 1)=0$ and $v(0)=0$, 
then the solution of \eqref{main00} is
\[
\Pi  f(\bx,t)=v(t)\prod_{j=1}^n w(x_j).
\]
If  $w \in C^N([p,q])$, we construct a Hestenes extension of $w(x)$ outside $[p,q]$ as
\begin{equation}\label{w}
\widetilde{w}(x)= \left\{ 
\begin{array}{cc}
\ds \sum_{s=1}^{N+1}c_s w(-\a_s(x-q)+q),& \ds q<x\leq q+\frac{q-p}{A} \\
w(x),& p\leq x \leq q\\
\ds\sum_{s=1}^{N+1}c_s w(-\a_s(x-p)+p),&  \ds p-\frac{q-p}{A}\leq x<p  \\
\end{array}
\right.
\end{equation}
where $\{a_1,...,a_{N+1}\}$  are different positive constants
$A=\max \, \a_s$,
and  ${\bf c}_N=\{c_1,...,c_{N+1}\}$ satisfy the system
\[
\sum_{s=1}^{N+1}c_s (-\a_s)^k=1,\quad k=0,...,N.
\]
Hence an extension of $f(\bx,t)$ with preserved smoothness is 
\[
\widetilde{f}(\bx,t)=v(t)\prod_{j=1}^n \widetilde w(x_j)\,.
\]
We compare the values of the exact and the approximate solution for
\eqref{main00}. In all the experiments the approximations have been
computed using \eqref{Sigmath3}-\eqref{Sigmath4} and the function
$\Psi_M$ in \eqref{PsiNew}. We choose the constants $\cD=\cD_0=4$ to
have the saturation error comparable with the double precision
rounding errors and the parameters in the quadrature rule
$\k=10^{-5}$, $R=3\cdot 10^6$, $a=1$.

In Tables \ref{one} and \ref{two} we report on the absolute errors and the approximation rates in the space dimensions $n=1,3,10,20,100,200$ for the solution of  \eqref{main00} with $v(t)=t$,  $w(x)=\cos^2(5 \pi x/2)$ and the Hestenes extension \eqref{w} with $\a_s=1/s$ (Table \ref{one}); $w(x)=\ee^{4 i x}(x^2-1)^2$ and $\widetilde{w}(x)=w(x)$  (Table \ref{two}). The results show that, for high dimensions, the second order fails but the forth and sixth order formulas approximate the exact solution  with the predicted approximation rates.

\begin{table}%[p]
%\begin{footnotesize}
\begin{center}
\caption{ Absolute errors and approximation rates
for the solution of \eqref{main00} with $f(\bx,t)$ in \eqref{fxt} where $w(x)=\cos^2(5 \pi x/2)$ and $v(t)=t$,  at the point $\bx=(0.1,0.4,...,0.4)$; $t=1$  using formula
\eqref{Sigmath3}-\eqref{Sigmath4} with \eqref{PsiNew} and the  Hestenes extension corresponding to $\a_s=1/s$.}\label{one}
\begin{tabular}{r|rr|cc|cc|cc} 
&&&\multicolumn{2}{c|}{$M=1$}&\multicolumn{2}{c|}{$M=2$}&\multicolumn{2}{c}{$M=3$}\\
&$h^{-1}$& $\t^{-1}$ &error  & rate         &         error     &         rate      &     error       &       rate     \\\midrule
 &$40$&$80$& 0.146E+00  &       &    0.326E-01  &      &   0.296E-02&             \\
$n=1$&$80$&$160$&  0.177E-01  &  3.04      & 0.106E-02 &   4.94   &      0.248E-04 &   6.89    \\
&$160$&$320$&  0.222E-02 &   2.99       &   0.313E-04  &  5.08    &    0.176E-06   & 7.13     \\
\midrule\midrule
 &$40$&$80$&  0.779E-01 &       &    0.135E-01  &     &    0.126E-02 &          \\
$n=3$&$80$&$160$&   0.194E-01   &  2.00        & 0.482E-03 &   4.80&    0.103E-04  &  6.93     \\
&$160$&$320$&   0.522E-02  &  1.89    &  0.240E-04  &  4.32     &    0.103E-04  &  6.93     \\
\midrule\midrule
 &$40$&$80$&   0.243E+00  &      &     0.236E-01  &   &   0.122E-02  &        \\
$n=10$&$80$&$160$&  0.789E-01  &  1.62      &  0.163E-02   & 3.86    &     0.208E-04  &  5.87      \\
&$160$&$320$&    0.212E-01  &  1.89  &     0.104E-03  &  3.97 &    0.356E-06   & 5.87       \\
\midrule\midrule
 &$40$&$80$&   0.378E+00  &      &     0.486E-01   &     &   0.258E-02  &           \\
$n=20$&$80$&$160$&   0.152E+00 &   1.31       & 0.343E-02   & 3.82      &  0.441E-04   & 5.87      \\
&$160$&$320$&    0.436E-01  &  1.80       &   0.219E-03   & 3.97   &    0.771E-06 &   5.84 \\
\midrule\midrule
 &$40$&$80$&   0.500E+00  &     &    0.207E+00  &     &    0.133E-01   &       \\
$n=100$&$80$&$160$&   0.424E+00  &  0.23       & 0.176E-01  &  3.55     &    0.230E-03   & 5.85      \\
&$160$&$320$&   0.189E+00  &  1.16    & 0.114E-02  &  3.95    &    0.402E-05  &  5.84       \\
\midrule\midrule
 &$40$&$80$&   0.500E+00  &     &   0.329E+00   &    &  0.264E-01  &     \\
$n=200$&$80$&$160$&   0.489E+00  &  0.03        & 0.348E-01  &  3.24   &      0.462E-03  &  5.84     \\
&$160$&$320$&  0.308E+00  &  0.66     &   0.229E-02 &   3.92     &  0.778E-05   & 5.89     \\
\bottomrule
\end{tabular}%\\[0.4mm]\\
\\
\bigskip
\caption{Absolute errors and approximation rates
for the solution of \eqref{main00} with $f(\bx,t)$ in \eqref{fxt} where $w(x)=\ee^{4 i x}(x^2-1)^2$ and $v(t)=t$,  
at the point $\bx=(0.1,0.1,...,0.1)$; $t=1$  using formula
\eqref{Sigmath3}-\eqref{Sigmath4} with \eqref{PsiNew} and the  extension $\widetilde{w}(x)=w(x)$
.}\label{two}
\begin{tabular}{r|rr|cc|cc|cc} 
&&&\multicolumn{2}{c|}{$M=1$}&\multicolumn{2}{c|}{$M=2$}&\multicolumn{2}{c}{$M=3$}\\
&$h^{-1}$& $\t^{-1}$ &error  & rate         &         error     &         rate      &     error       &       rate     \\\midrule
&$20 $&$40$& 0.638E-01  &         &   0.153E-02&     &    0.724E-04  &   \\
$n=1$&$40$&$80$&  0.162E-01    &1.98      &  0.986E-04  &  3.96     &      0.122E-05 &   5.89   \\
&$80 $&$160$&  0.407E-02   & 1.99         &  0.621E-05  &  3.99     &      0.199E-07 &   5.94   \\
\midrule\midrule
&$20 $&$40$&0.133E+00   &         &0.550E-02   &     &   0.168E-03   &   \\
$n=3$&$40$&$80$&   0.354E-01 &    1.96    &  0.361E-03  &  3.93     &   0.277E-05   & 5.92   \\
&$80 $&$160$&  0.899E-02  &  1.97        &   0.228E-04  &  3.98     &      0.439E-07 &   5.98   \\
\midrule\midrule
&$20 $&$40$&  0.321E+00 &         & 0.161E-01  &     &    0.512E-03  &   \\
$n=10$&$40$&$80$&  0.968E-01  &  1.73       &   0.106E-02  &  3.92     &     0.843E-05   & 5.92 \\
&$80 $&$160$& 0.254E-01  &  1.92        &  0.672E-04 &   3.98     &     0.134E-06  &  5.98   \\
\midrule\midrule
&$20 $&$40$& 0.423E+00  &         & 0.260E-01  &     & 0.837E-03     &   \\
$n=20$&$40$&$80$&   0.149E+00    & 1.50         &  0.174E-02 &   3.91    &      0.138E-04  &  5.92  \\
&$80 $&$160$&   0.409E-01   & 1.86        &   0.110E-03  &  3.98     &    0.219E-06  &  5.98   \\
\midrule\midrule
&$20 $&$40$&  0.133E+00 &         & 0.242E-01  &     &     0.836E-03 &   \\
$n=100$&$40$&$80$&   0.964E-01   & 0.46        &   0.173E-02   & 3.80   &      0.138E-04   & 5.92   \\
&$80 $&$160$&   0.363E-01  &  1.41         &   0.110E-03  &  3.97    &   0.219E-06  &  5.98  \\
\midrule\midrule
&$20 $&$40$&  0.180E-01  &         &0.590E-02   &     &   0.223E-03   &   \\
$n=200$&$40$&$80$&   0.166E-01   & 0.11       &  0.461E-03 &   3.68    &      0.370E-05 &    5.91 \\
&$80 $&$160$&    0.843E-02  &  0.97       &   0.295E-04 &   3.97     &    0.587E-07   & 5.98   \\
\bottomrule
\end{tabular}%\\[0.4mm]
\end{center}
%\end{footnotesize}
\end{table}

\subsection{Initial value problem} \label{subsec52}

Consider the initial value problem 
\begin{equation}\label{ex1}
i \frac{\partial u}{\partial t} +\D_\bx u=0,\quad u(\bx,0)=g(\bx)=\prod_{j=1}^n w(x_j),\quad w(x_j)=0\quad \hbox{ if} \quad x_j \not\in [-1,1].
\end{equation}
Thus $\supp g\subset [-1,1]^n$. Denote by $\widetilde{w}$ the extension of $w$ outside $[-1,1]$ with preserved smoothness. An approximate solution of \eqref{ex1} is given by
\begin{equation}\label{app4}
\begin{split}
 u_h(\bx,t)&=\frac{1}{\cD^{n/2}}\prod_{j=1}^n \sum_{h m \in I }\widetilde{w}(h m) \\
 &\times\Big(
 \Psi_M\big(\frac{x_j-hm}{h\sqrt{\cD}},\frac{4t}{h^2\cD},\frac{-1-hm}{h\sqrt{\cD}}\big)-
\Psi_M\big(\frac{x_j-hm}{h\sqrt{\cD}},\frac{4t}{h^2\cD},\frac{1-hm}{h\sqrt{\cD}}\big)
 \Big)
 \end{split}
  \end{equation}
with $I=(-1-r \sqrt{\cD},1+r \sqrt{\cD})$.

\begin{figure}
\centering
 \includegraphics[scale=0.6]{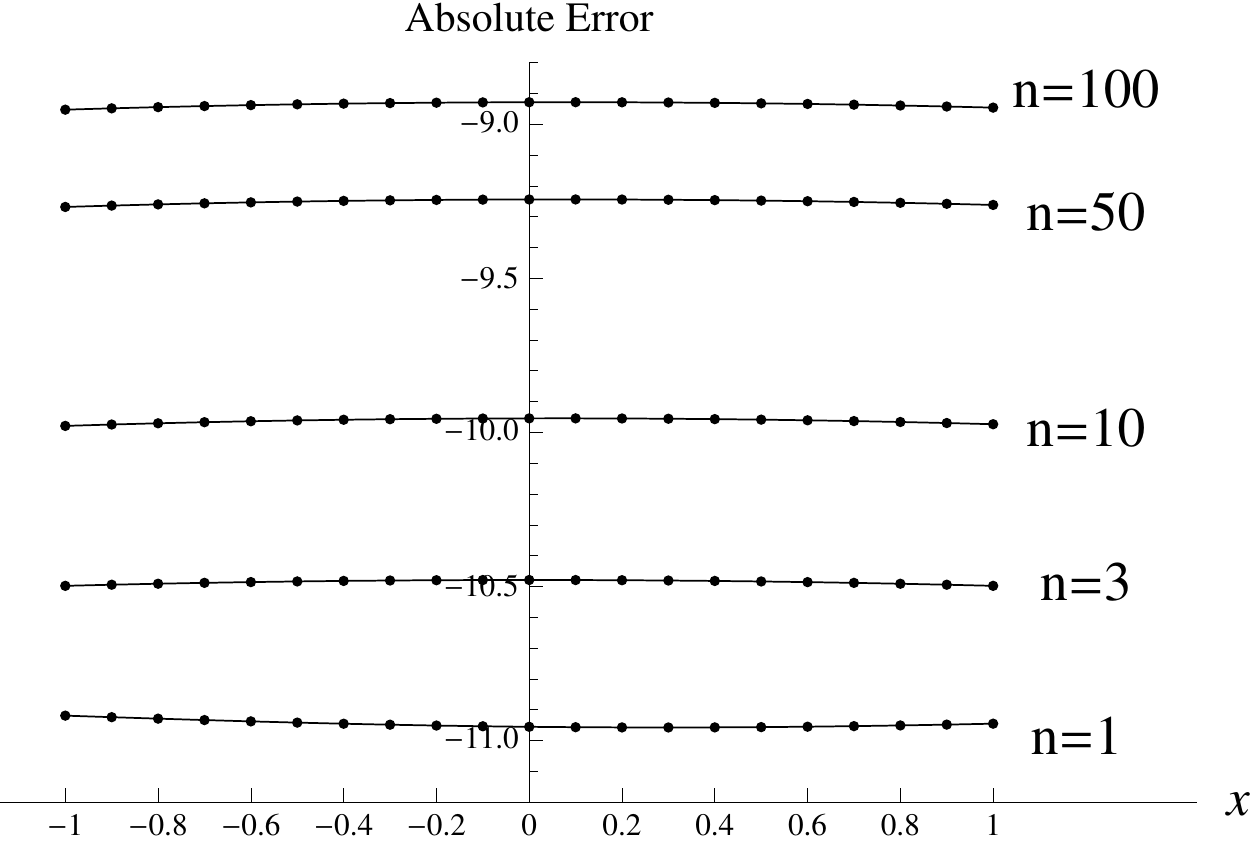}
\caption{ {\footnotesize Absolute errors, using $\log_{10}$ scale on the vertical
axes, for the solution of \eqref{ex1} with $w(x)=\ee^{(x+a)^2}$, $a=0.32612$, the Hestenes  extension
corresponding to $\a_s=1/2^s$, using \eqref{app4} with $h=1/160$, $\cD=4$, $\bx=(x,0.1,...,0.1)$, $t=1$.}} \label{fig3}
\end{figure}

In this part we provide results of some experiments which show accuracy and numerical convergence orders.
We assume $w(x)=\ee^{(x+a)^2}$ which gives the exact solution of \eqref{ex1}  
\[
u(\bx,t)=\prod_{j=1}^n \frac{i \ee^{\frac{(a+x_j)^2}{1-4 i t}}}{{2 \sqrt{4 i t-1}}}
 \left({\erfc}\left(\frac{4 i (a+1) t+x_j-1}{2 \sqrt{t} \sqrt{4
   t+i}}\right)-{\erfc}\left(\frac{4 i (a-1) t+x_j+1}{2 \sqrt{t} \sqrt{4 t+i}}\right)\right)\,
\]
and we compare the calculated solution $u_h$ with the exact solution
$u$.  In our experiments we choose $a=0.32612$.  In Figure \ref{fig3}
we report on the absolute error at some grid points in dimensions
$n=1,3,10,50,100$. The approximations have been computed with $\cD=4$,
$M=3$ and $h=1/160$ in \eqref{app4}, and the Hestenes extension with
$\a_s=1/2^s$.  If $g$ allows the representation \eqref{sep} that is
$g$ has rank $P$ then, denoting by $\e^{(p)}_j$ the 1-dimensional
error for each function $g_j^{(p)}$, then the total error
$\e_n=\cO\left(\sum_{p=1}^P \sum_{j=1}^n \e_j^{(p)}\right)$. Results
in Figure \ref{fig3} confirm that, for $P=1$ and $g$ in \eqref{ex1},
the $n-$dimensional error $\e_n=\cO(n \e_1)$.

In Table \ref{three} we show that formula \eqref{app4} approximates
the exact solution with the predicted approximate orders $N=2,4,6$ in
the space dimensions $n=1,3,10,20,\\100,200$.

\begin{table}
%\begin{footnotesize}
\begin{center}
\caption{ Absolute errors and approximation rates
for the solution of \eqref{ex1} with $w(x)=\ee^{(x+a)^2}$,  $a=0.32612$, at the point $\bx=(0.2,0.1,...,0.1)$; $t=1$  using formula
\eqref{app4}  and the  Hestenes extension corresponding to $\a_s=1/s$
.}\label{three}
\begin{tabular}{r|r|cc|cc|cc} 
&&\multicolumn{2}{c|}{$M=1$}&\multicolumn{2}{c|}{$M=2$}&\multicolumn{2}{c}{$M=3$}\\
&$h^{-1}$&error  & rate         &         error     &         rate      &     error       &       rate     \\\midrule

 &$40$& 3.069E-03     &       &   1.178E-05    &      &4.522E-08    &             \\
 $n=1$&$80$&   7.693E-04  & 1.99      &      7.438E-07  & 3.98      &   7.206E-10   &5.97           \\
 &$160$&      1.924E-04   & 1.99     &   4.661E-08  & 3.99  &   1.151E-11 &  5.96            \\
 &$320$&  4.812E-05 &  1.99      &     2.915E-09 &  3.99    &   2.158E-13  & 5.73             \\
 
\midrule\midrule

 &$40$&   9.246E-03   &      & 3.538E-05      &    & 1.357E-07     &         \\
$n=3$ &$80$&  2.312E-03  & 1.99    &     2.233E-06 &  3.98      &  2.163E-09  & 5.97            \\
 &$160$&  5.781E-04  & 1.99     &      1.399E-07  & 3.99      &    3.457E-11 &  5.96             \\
 &$320$&  1.445E-04  & 1.99       &      8.754E-09  & 3.99      &    6.455E-13 &   5.74             \\
\midrule\midrule

 &$40$&   1.292E-01  &     &1.184E-04      &      & 4.546E-07   &             \\
$n=10$ &$80$&  7.764E-03   & 2.01       &     7.479E-06   &3.98      &    7.246E-09  & 5.97             \\
 &$160$&  1.937E-03  & 2.00      &      4.687E-07  & 3.99      &    1.157E-10 &  5.96          \\
 &$320$&   4.840E-04  & 2.00      &      2.931E-08  & 3.99    &   2.159E-12  & 5.74            \\
 \midrule\midrule
 
  &$40$&  6.400E-02     &      & 2.385E-04     &    &9.155E-07 &             \\
$n=20$ &$80$&  1.569E-02 &   2.02      &      1.505E-05  & 3.98      &    1.458E-08  & 5.97            \\
 &$160$&  3.904E-03 &  2.00       &     9.437E-07  & 3.99     &    2.331E-10 &  5.96           \\
 &$320$&  9.749E-04&   2.00      &      5.902E-08   & 3.99     &  4.347E-12 &  5.74           \\
 \midrule\midrule

 &$40$& 3.832E-01  &       & 1.258E-03      &     &    4.831E-06 &            \\
$n=100$ &$80$&  8.542E-02  & 2.16       &    7.947E-05 &  3.98     &   7.700E-08   &5.97            \\
 &$160$&  2.076E-02  & 2.04       &      4.980E-06   &3.99      &   1.230E-09   &5.96           \\
 &$320$& 5.155E-03 &  2.01       &     3.115E-07 &  3.99      &  2.294E-11&   5.74           \\ 
 \midrule\midrule

 &$40$&  9.669E-01 &       &  2.691E-03    &      &1.033E-05 &             \\
$n=200$ &$80$&   1.901E-01   &     2.34     &      1.700E-04 &  3.98      &    1.647E-07 &  5.97            \\
 &$160$&   4.486E-02  & 2.08       &      1.065E-05 &  3.99    &    2.632E-09&   5.96         \\
 &$320$&1.105E-02  & 2.02      &      6.665E-07  & 3.99      &    4.908E-11   &5.74           
   \\\bottomrule
\end{tabular}%\\[0.4mm]
\end{center}
%\end{footnotesize}
\end{table}

\begin{figure}[p]
\includegraphics[scale=0.5]{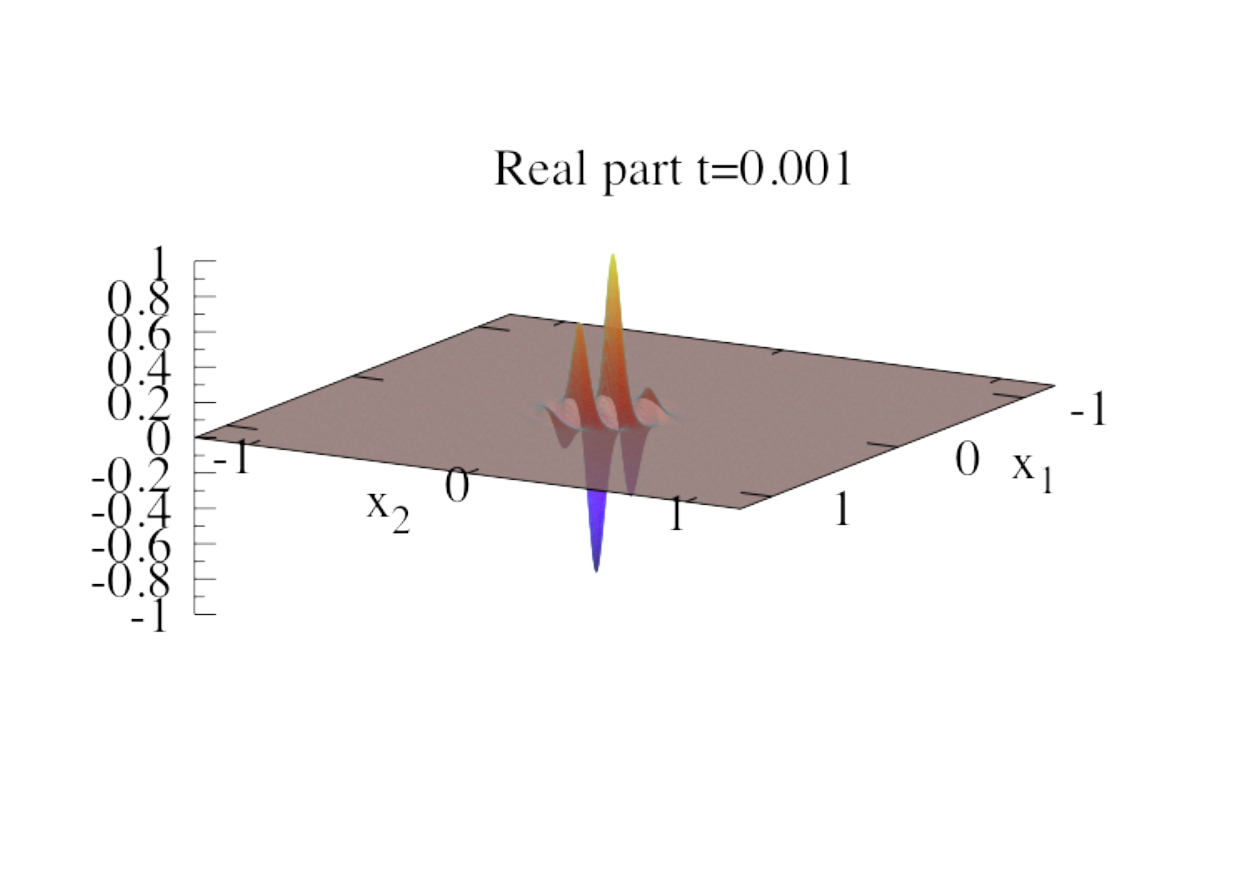}   
\includegraphics[scale=0.5]{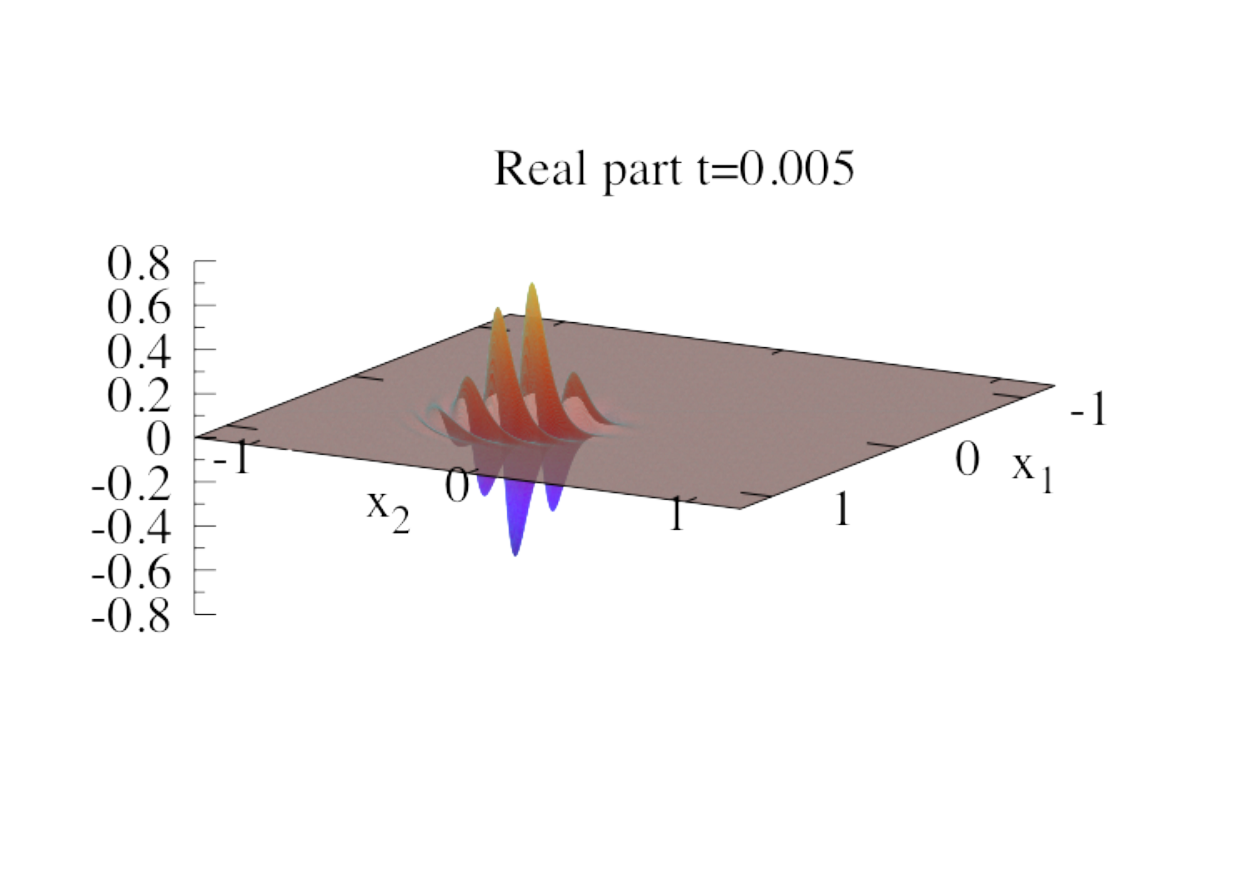}\\   [-1cm]
\includegraphics[scale=0.5]{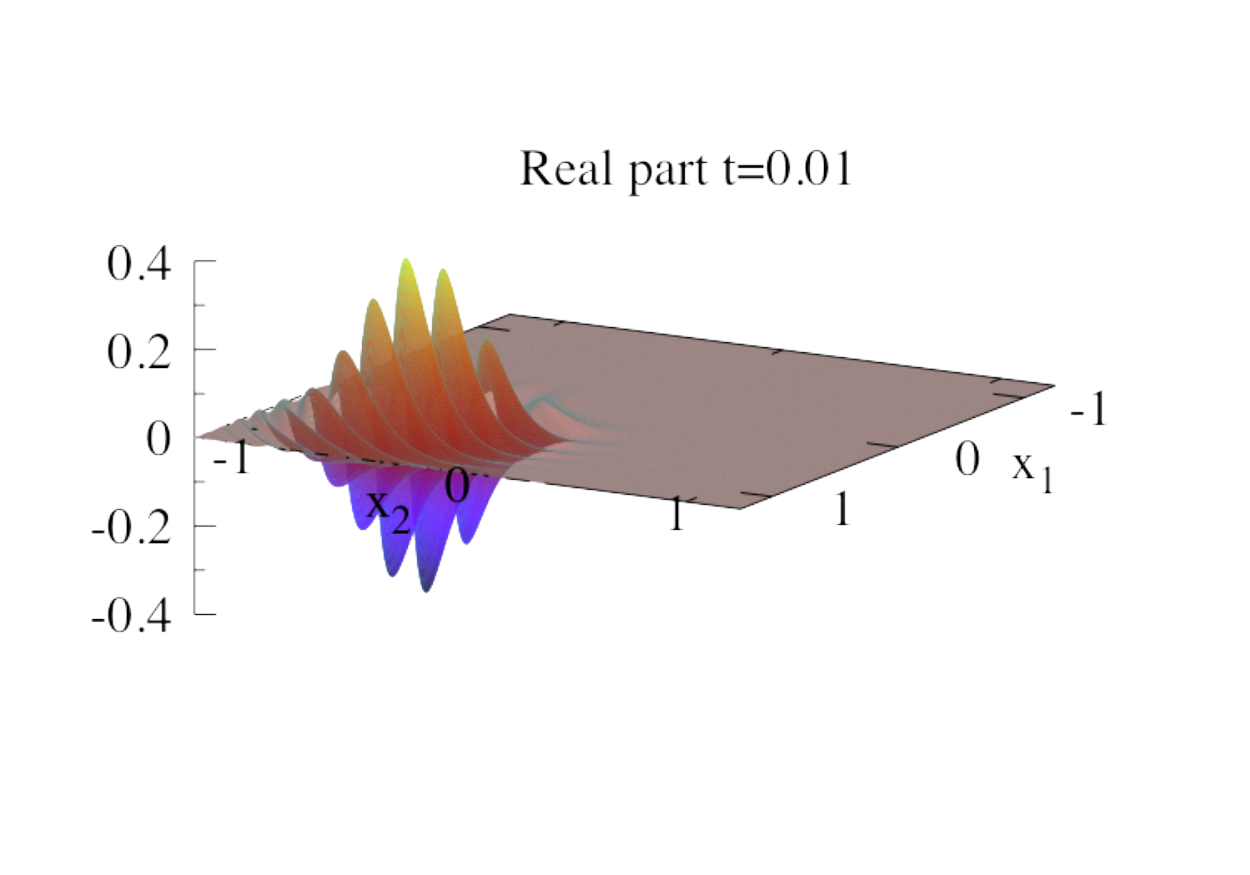}  
\includegraphics[scale=0.5]{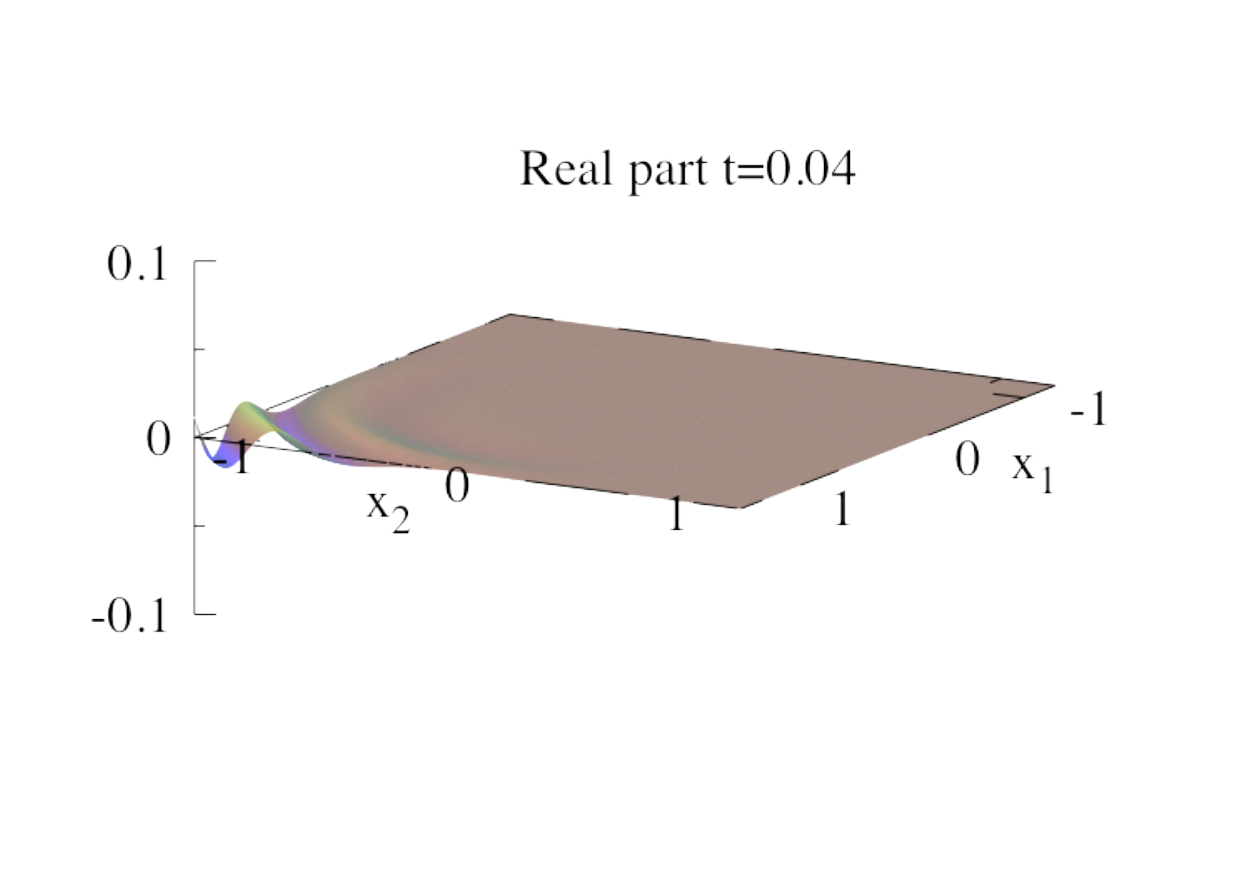} 
\vspace{-1.2cm}
\caption{Real part of $u$ when $w_j(x)=\ee^{i c_jx} \ee^{-60x^2}$, $j=1,2$, $c_1=30$, $c_2=-30$, $N=6$, $h=0.005$\,.}\label{figre30}
\includegraphics[scale=0.5]{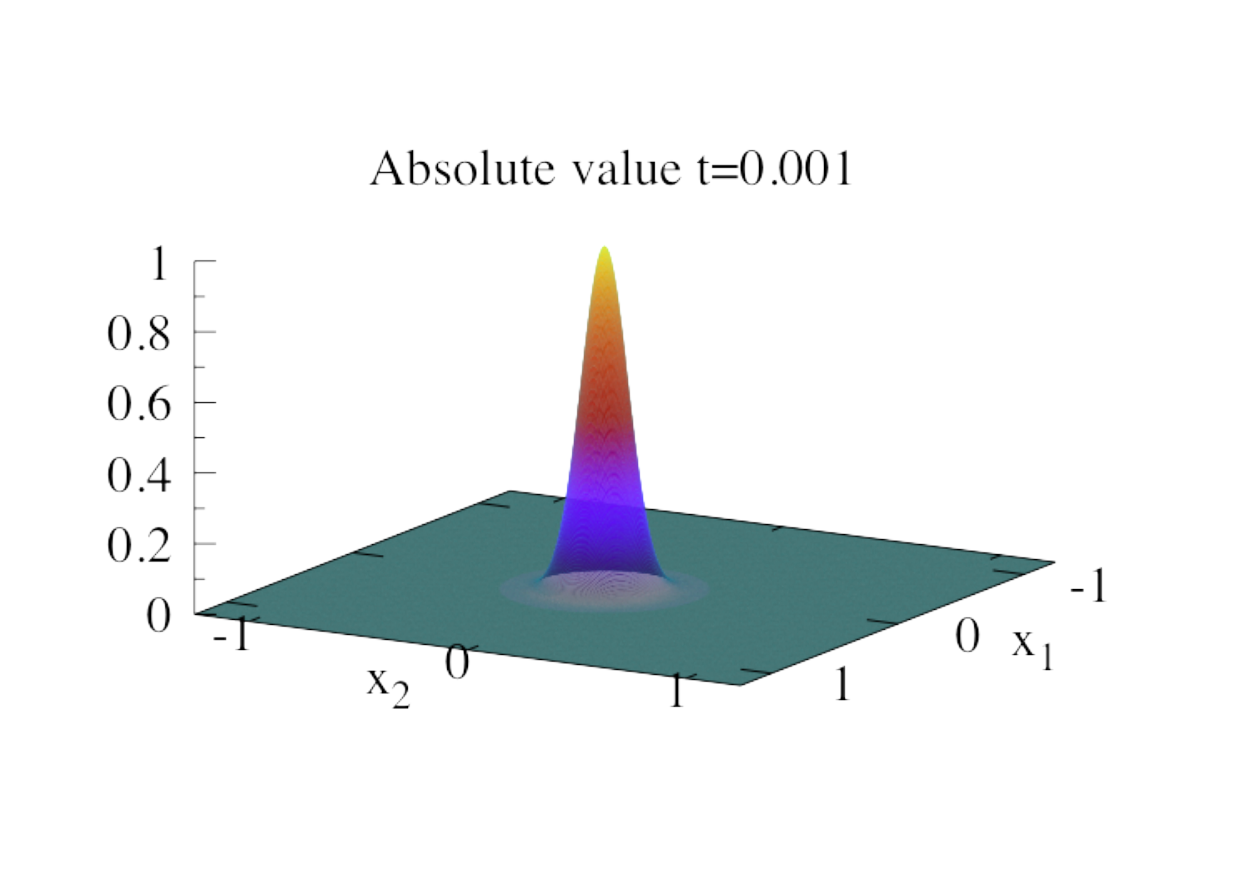}  
\includegraphics[scale=0.5]{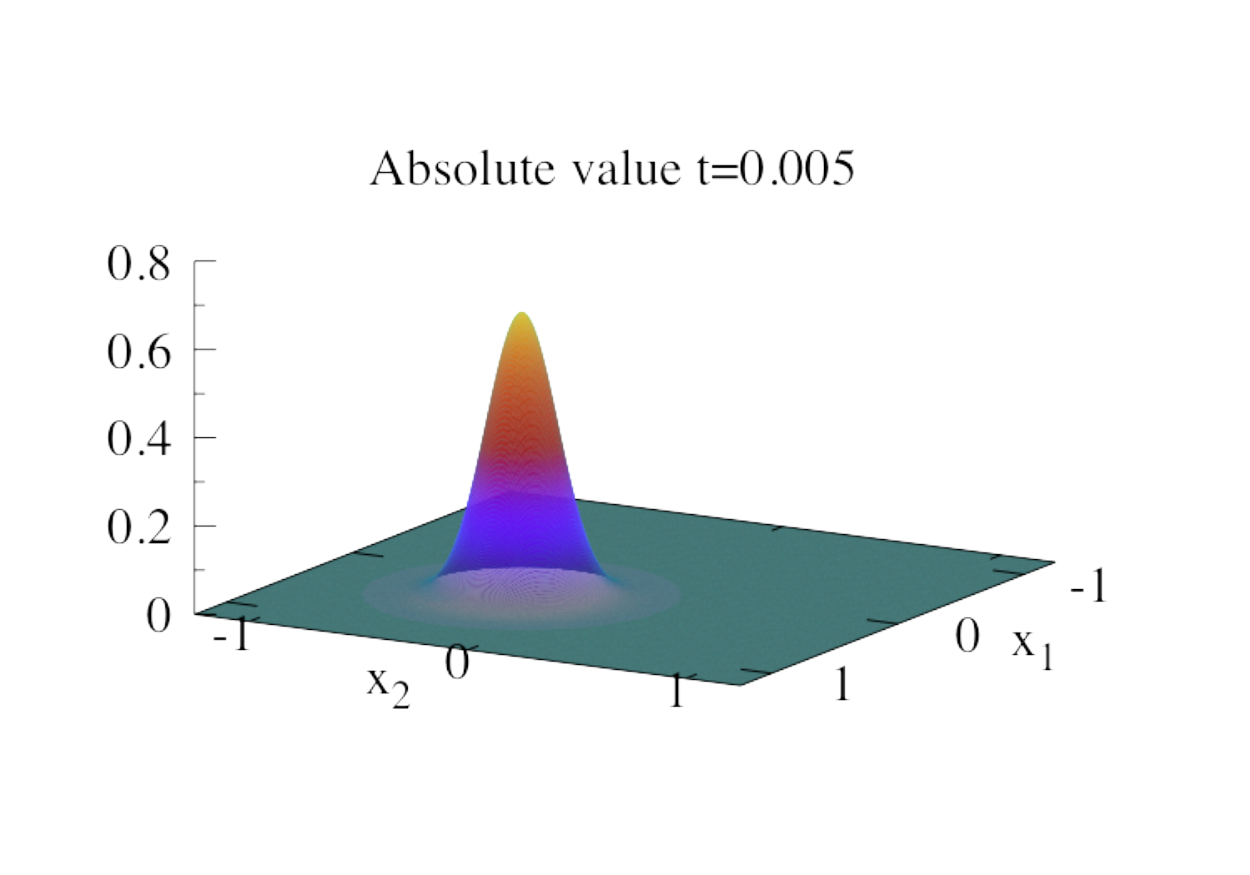}\\ [-1cm]
\includegraphics[scale=0.5]{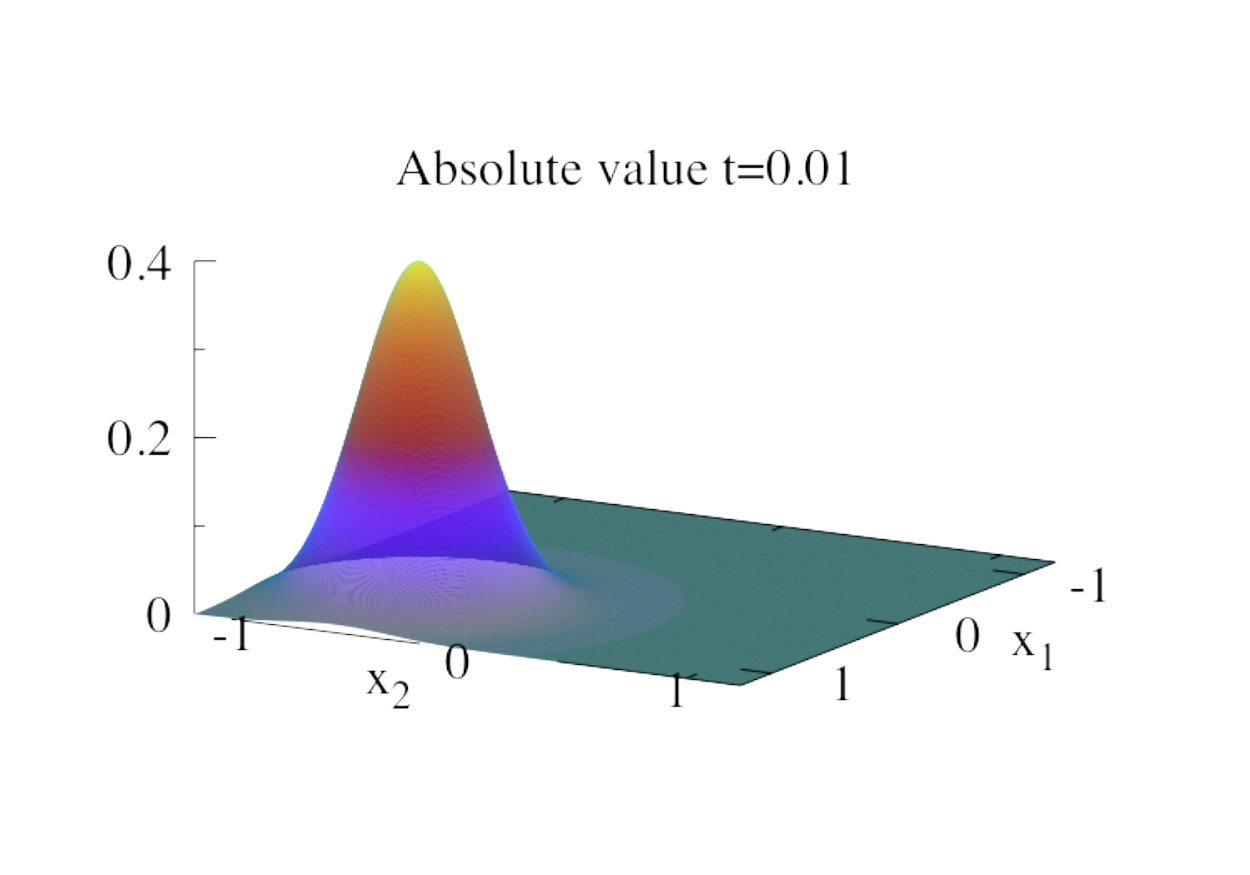}  
\includegraphics[scale=0.5]{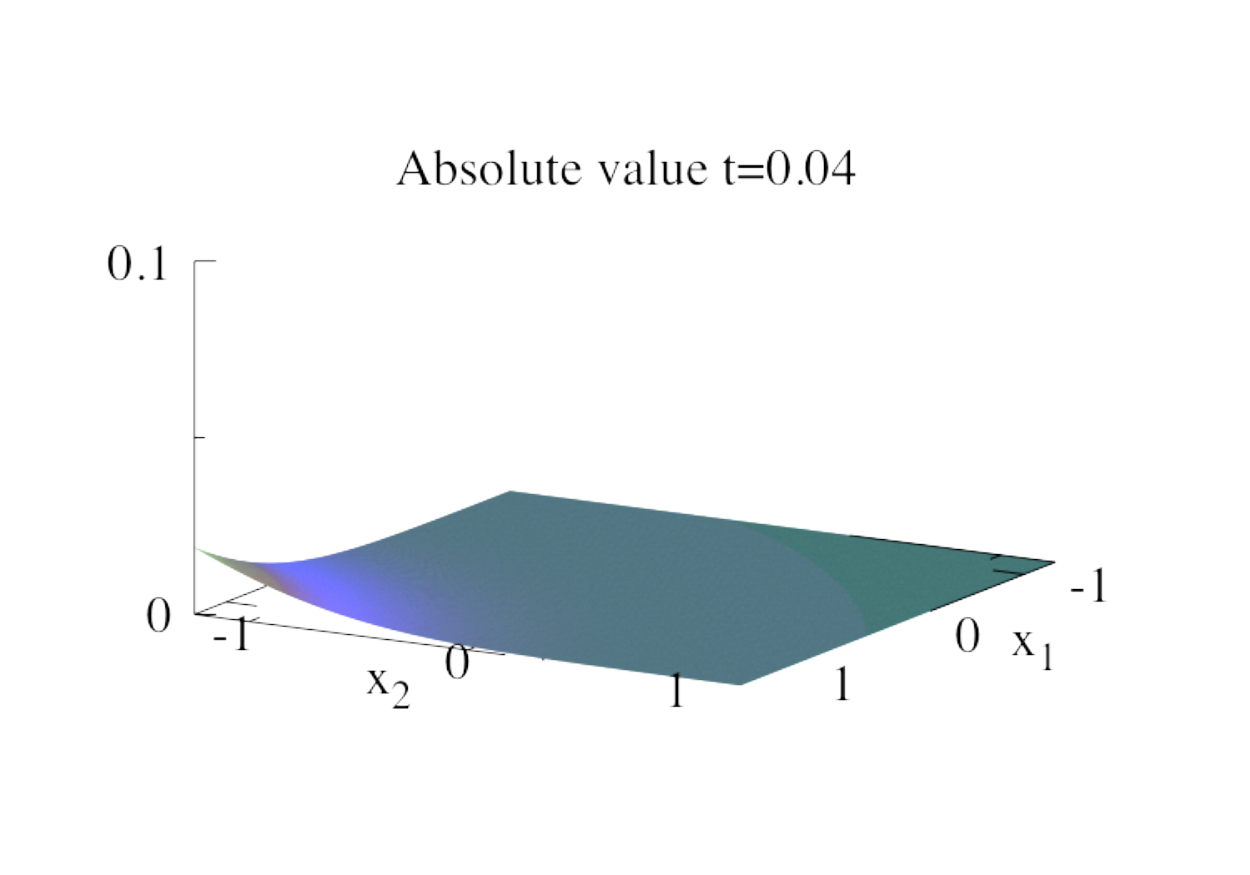} 
\vspace{-1.2cm}
\caption{Absolute value of $u$ when $w_j(x)=\ee^{i c_jx} \ee^{-60x^2}$, $j=1,2$, $c_1=30$, $c_2=-30$, $N=6$, $h=0.005$\,.}\label{figab30}
\end{figure}

\begin{figure}[p]
\includegraphics[scale=0.5]{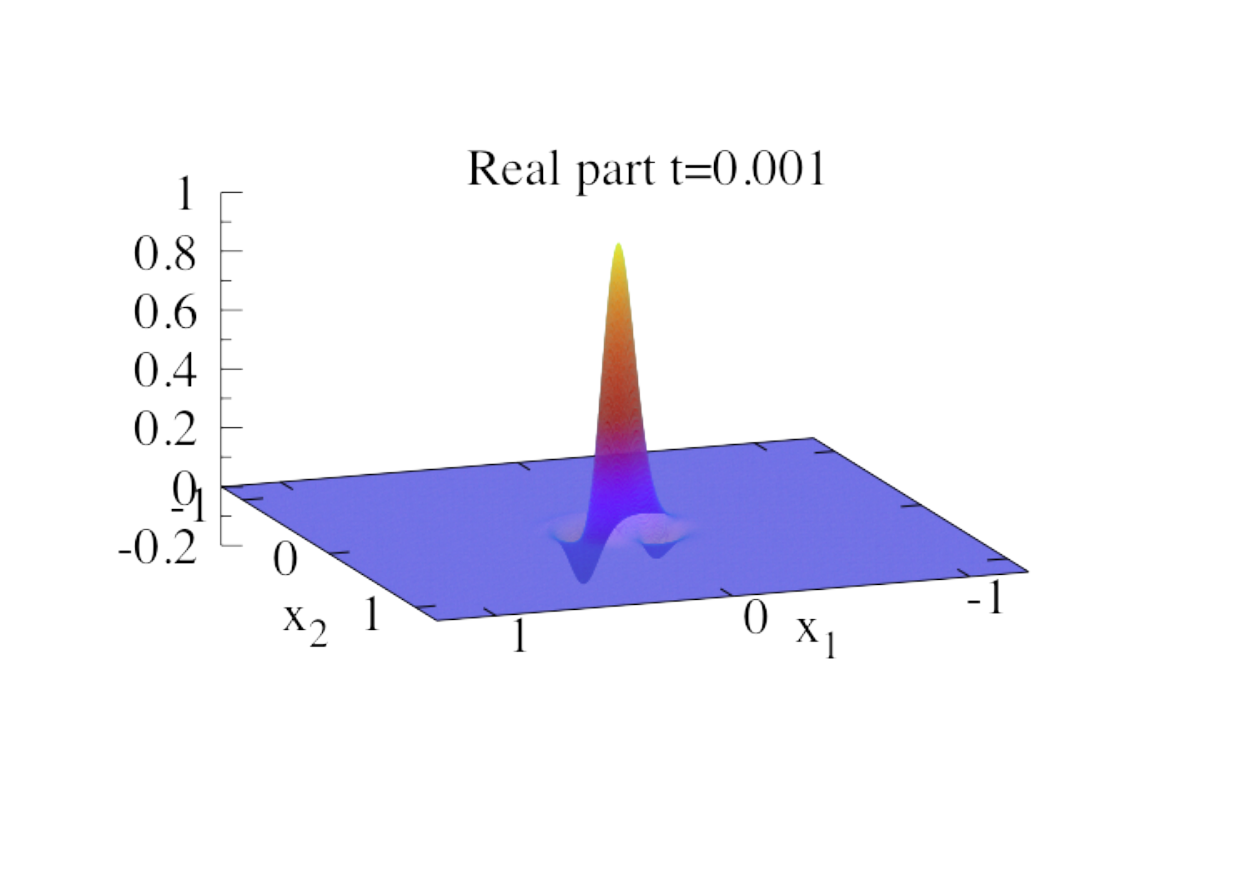}
\includegraphics[scale=0.5]{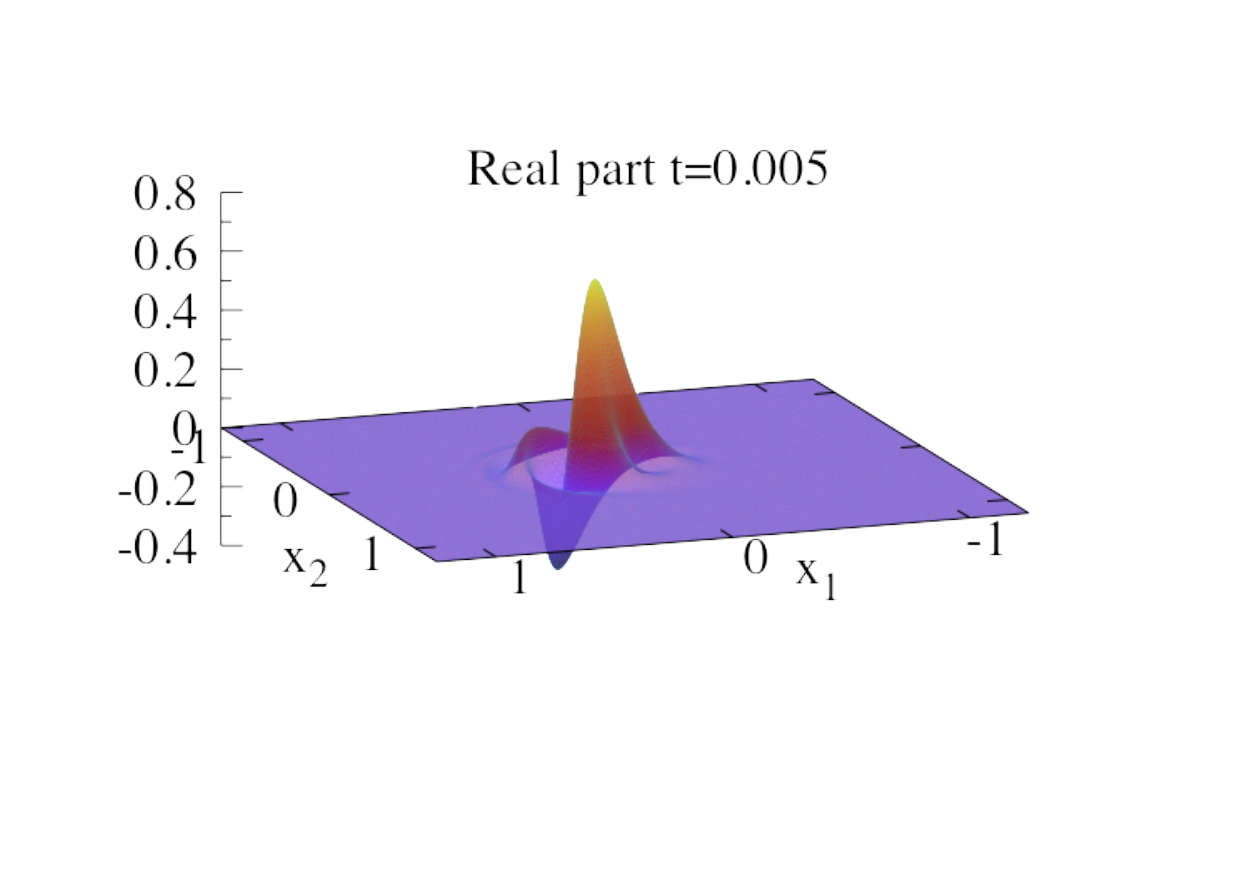}\\[-1cm]
\includegraphics[scale=0.5]{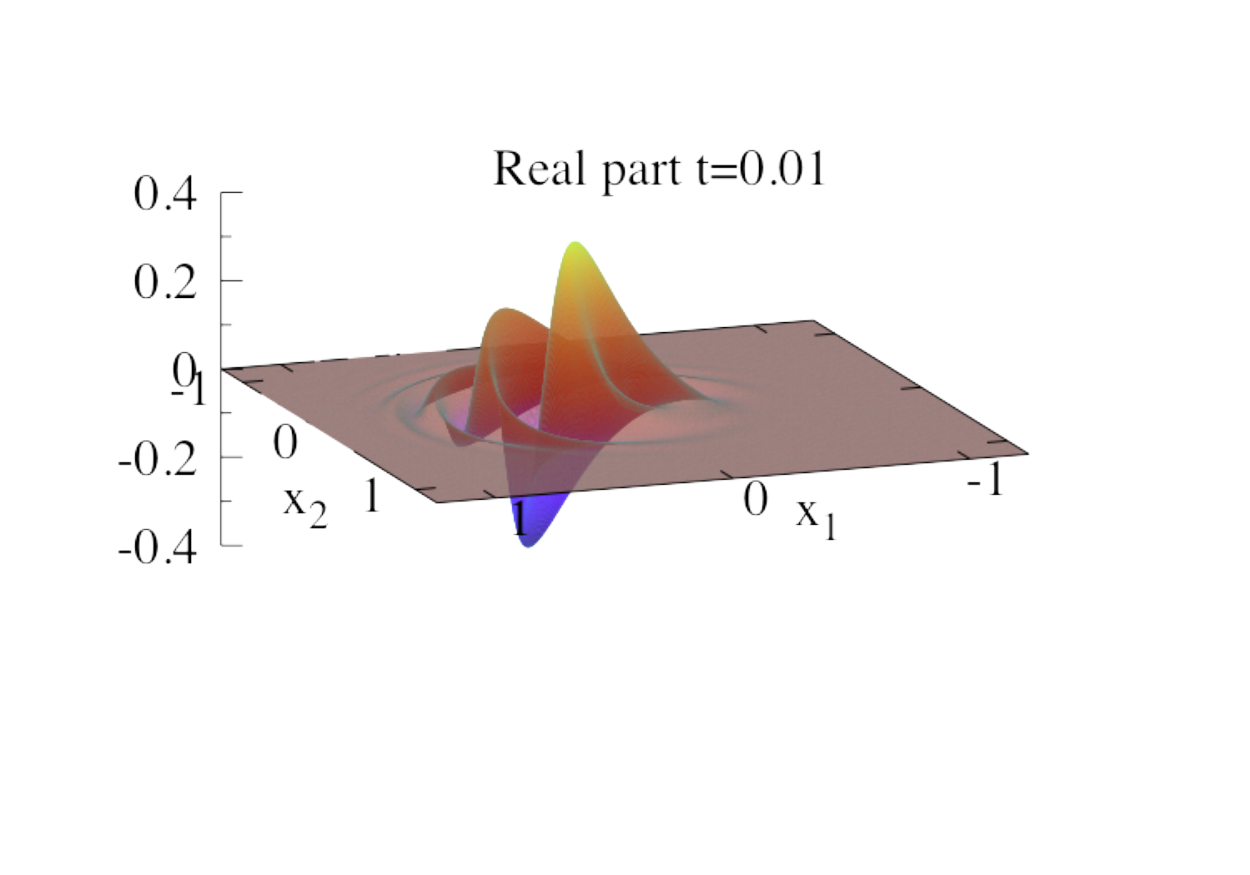}
\includegraphics[scale=0.5]{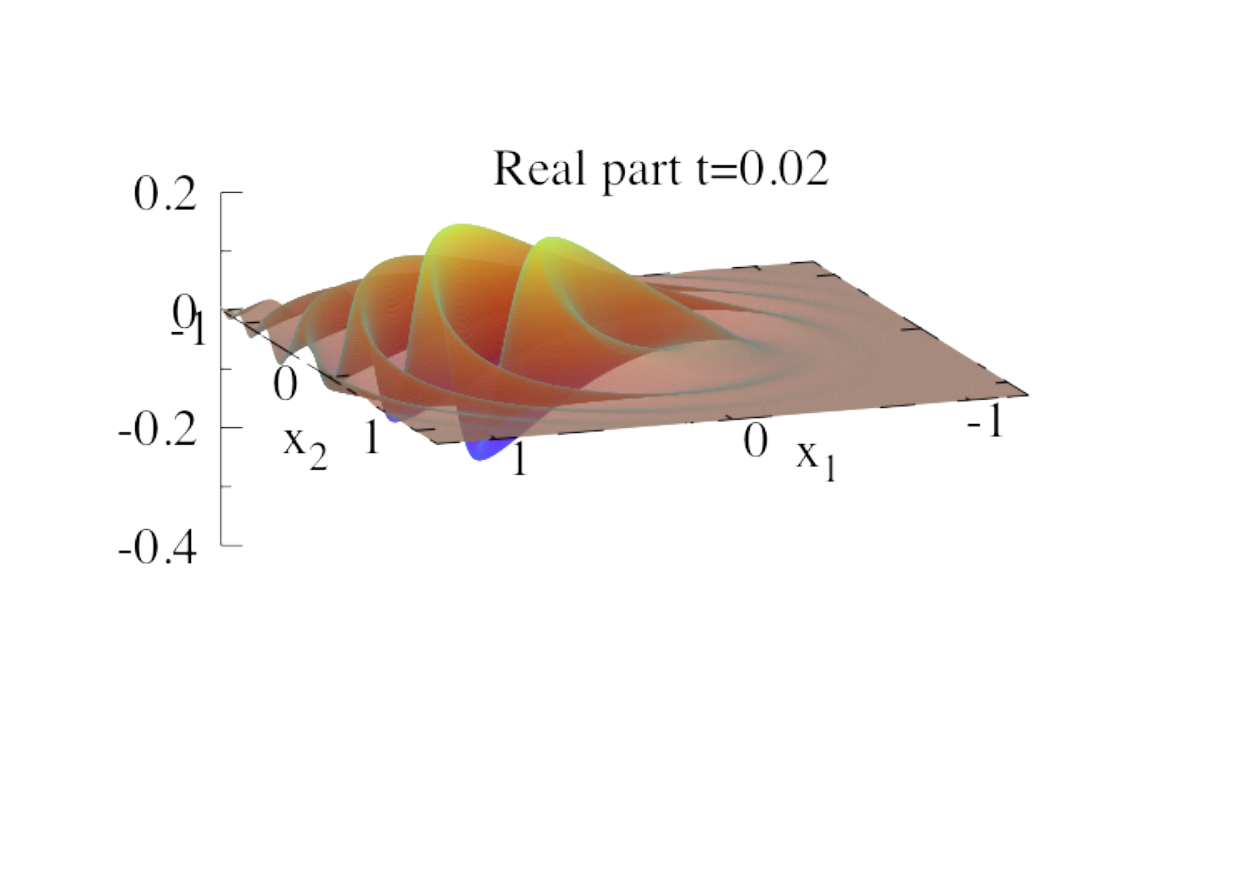}
\vspace{-1.2cm}
\caption{Real part of $u$ when $w_j(x)=\ee^{i c_jx} \ee^{-60x^2}$, $j=1,2$, $c_1=10$, $c_2=-10$, $N=6$, $h=0.005$\,.}\label{figre10}
\includegraphics[scale=0.5]{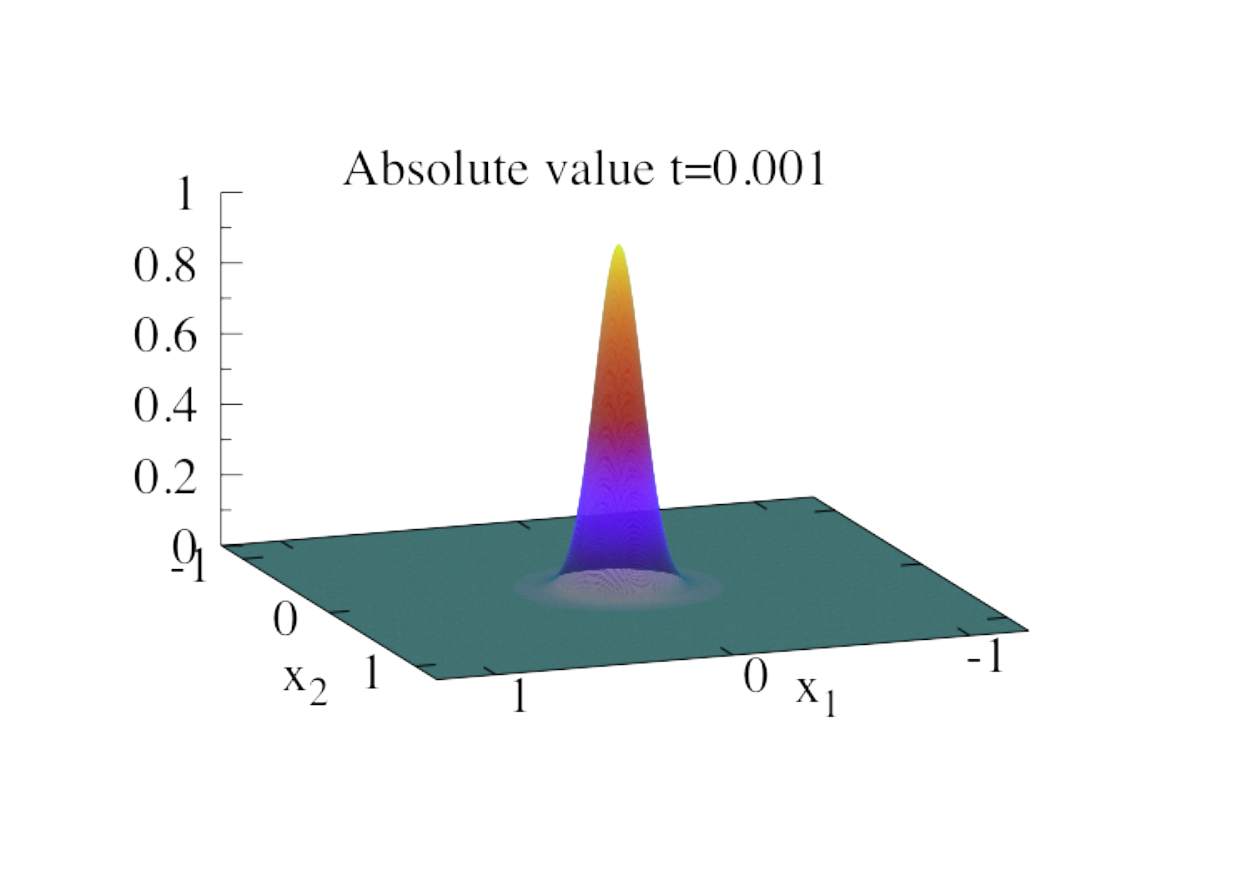}  
\includegraphics[scale=0.5]{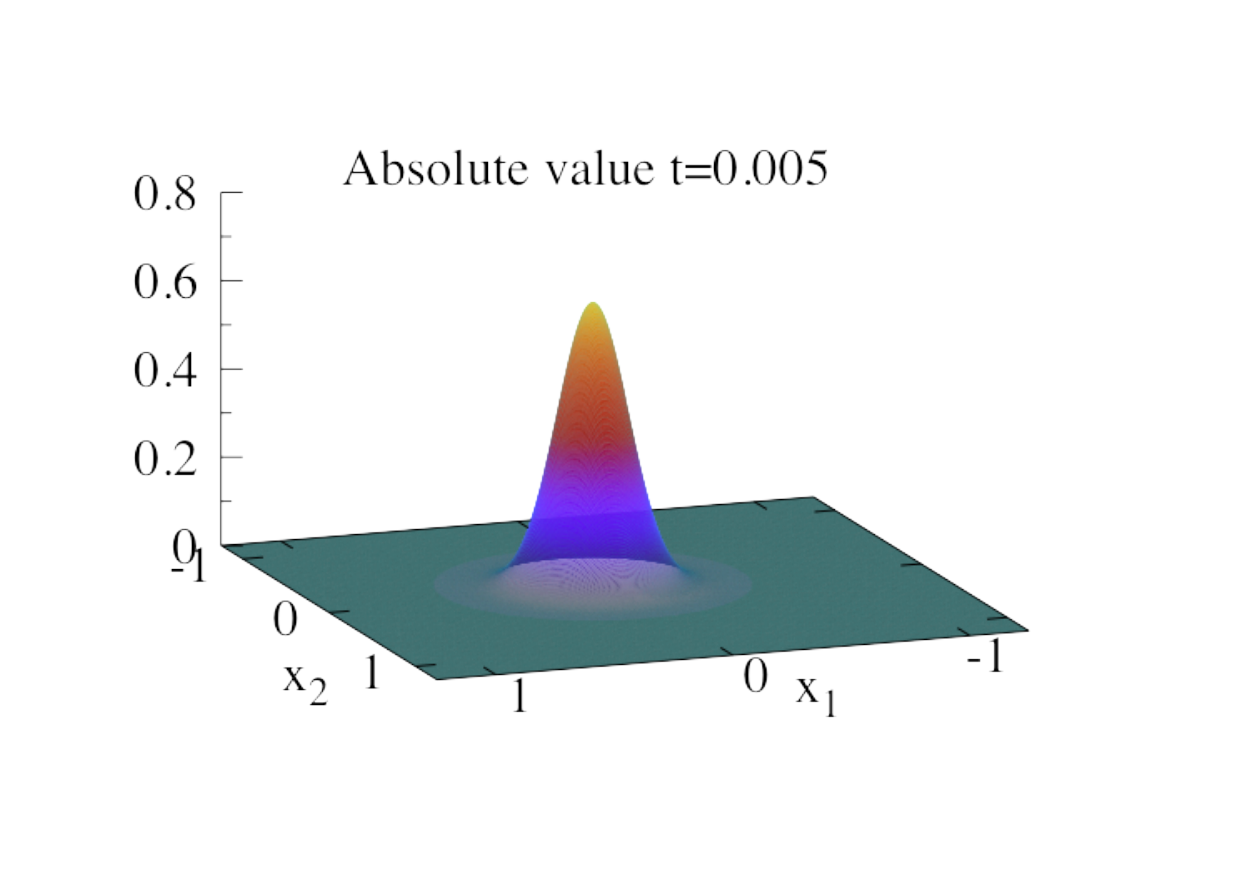} \\[-1cm]
\includegraphics[scale=0.5]{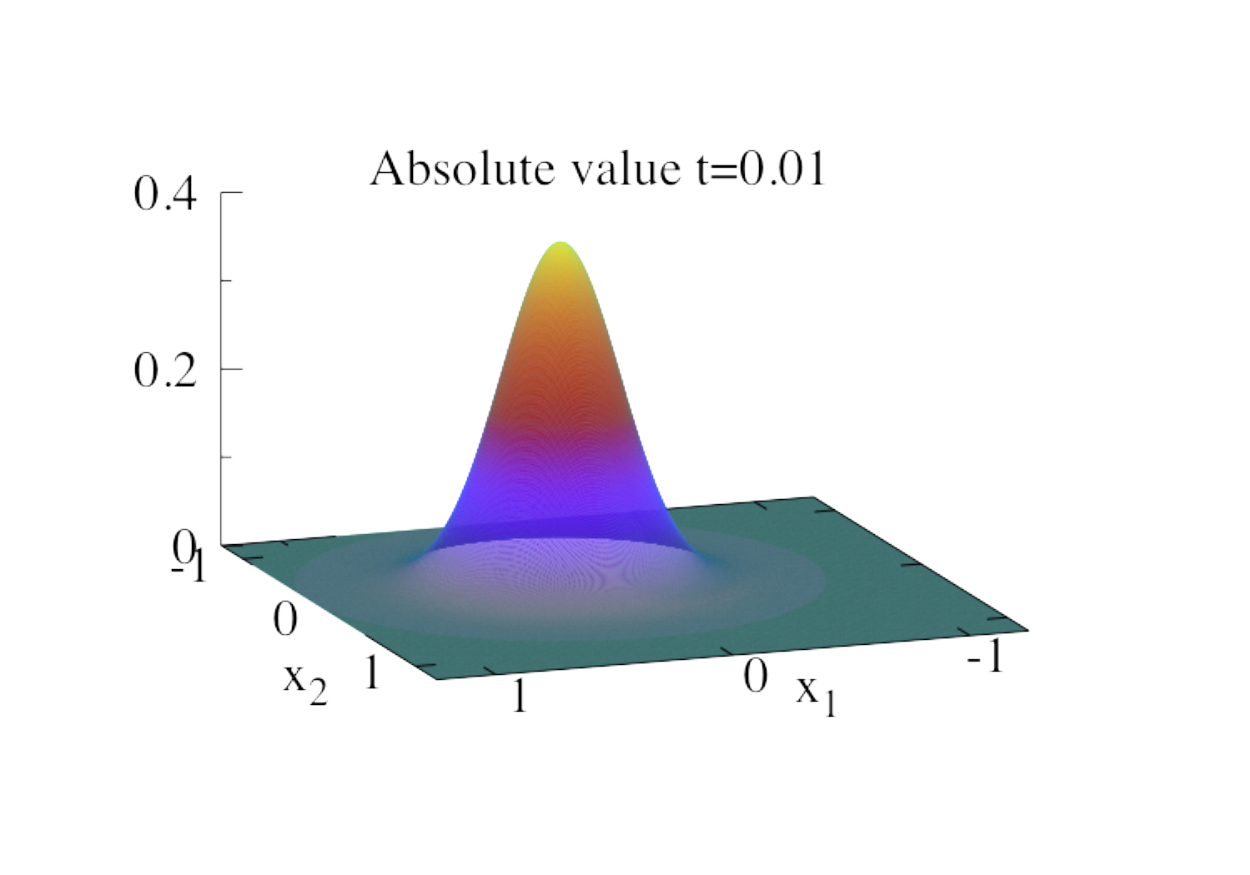}  
\includegraphics[scale=0.5]{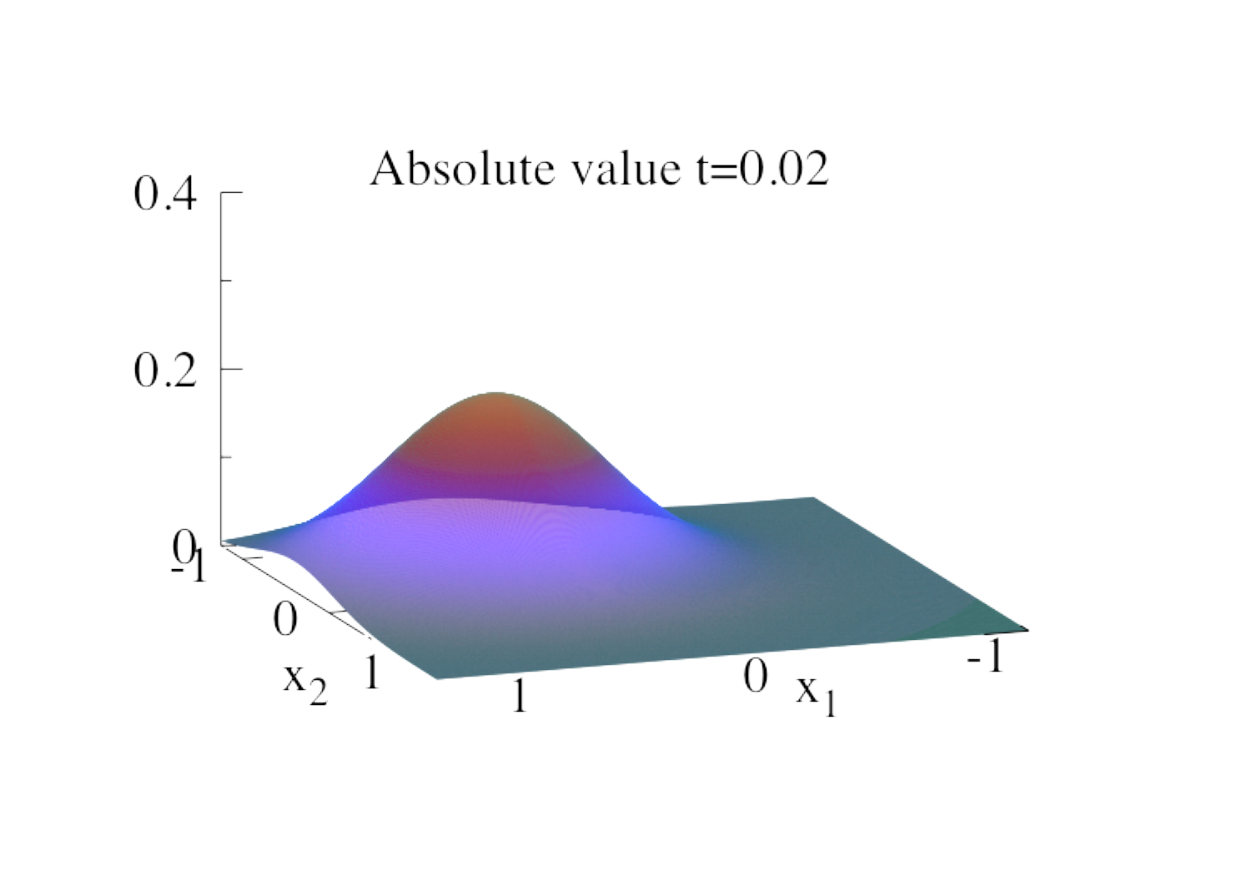}  
\vspace{-1.2cm}
\caption{Absolute value of $u$ when $w_j(x)=\ee^{i c_jx} \ee^{-60x^2}$, $j=1,2$, $c_1=10$, $c_2=-10$, $N=6$, $h=0.005$}\label{figab10}
\end{figure}

\begin{figure}[p]
\includegraphics[scale=0.5]{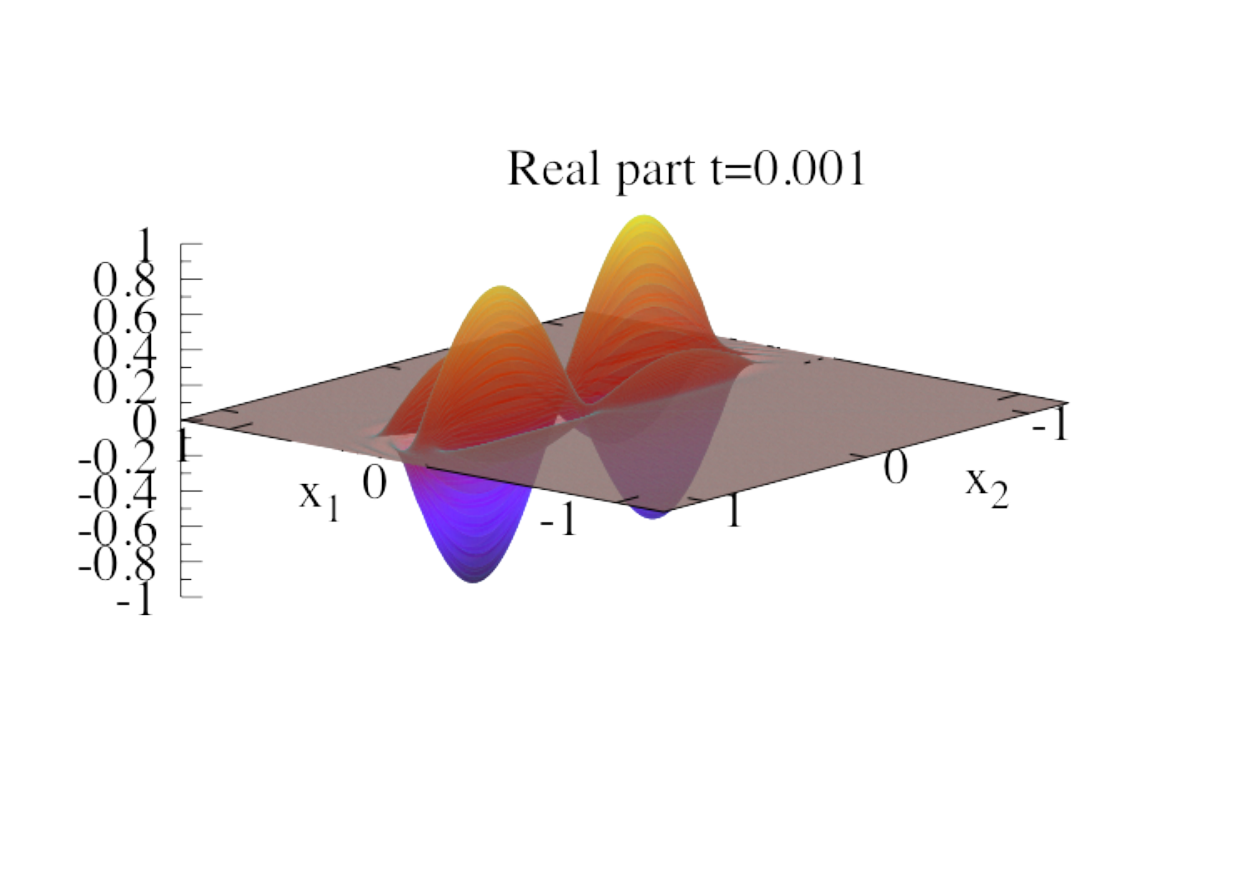} 
\includegraphics[scale=0.5]{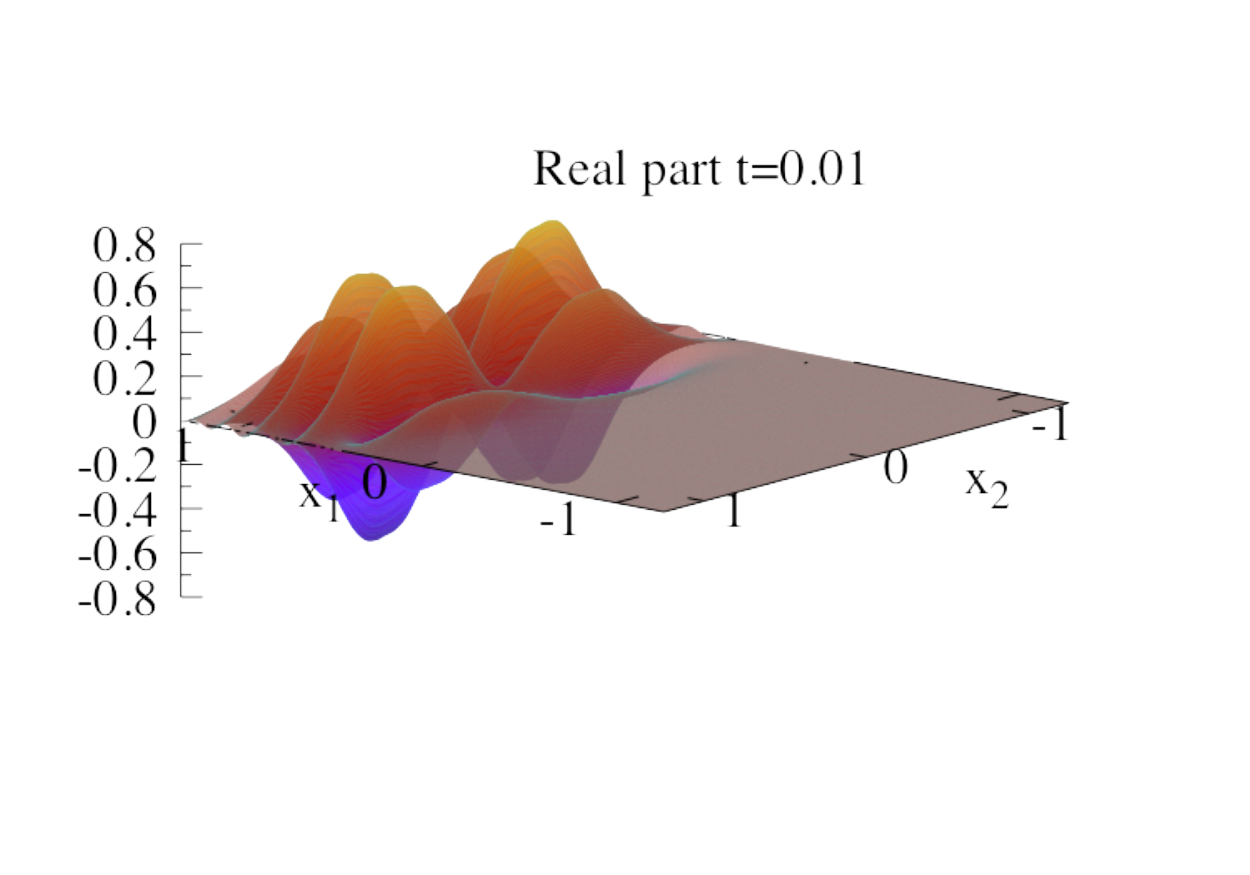}\\ [-1cm]
\includegraphics[scale=0.5]{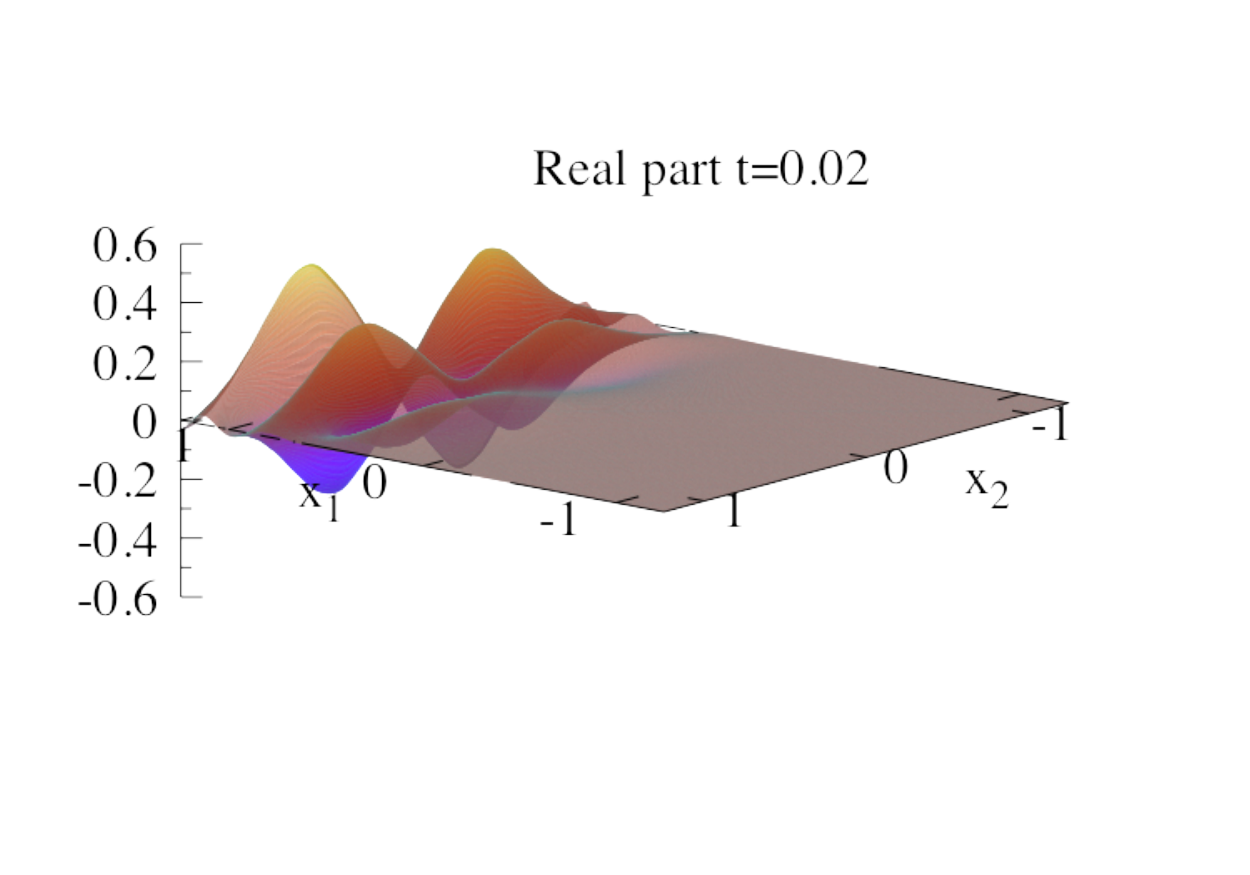} 
\includegraphics[scale=0.5]{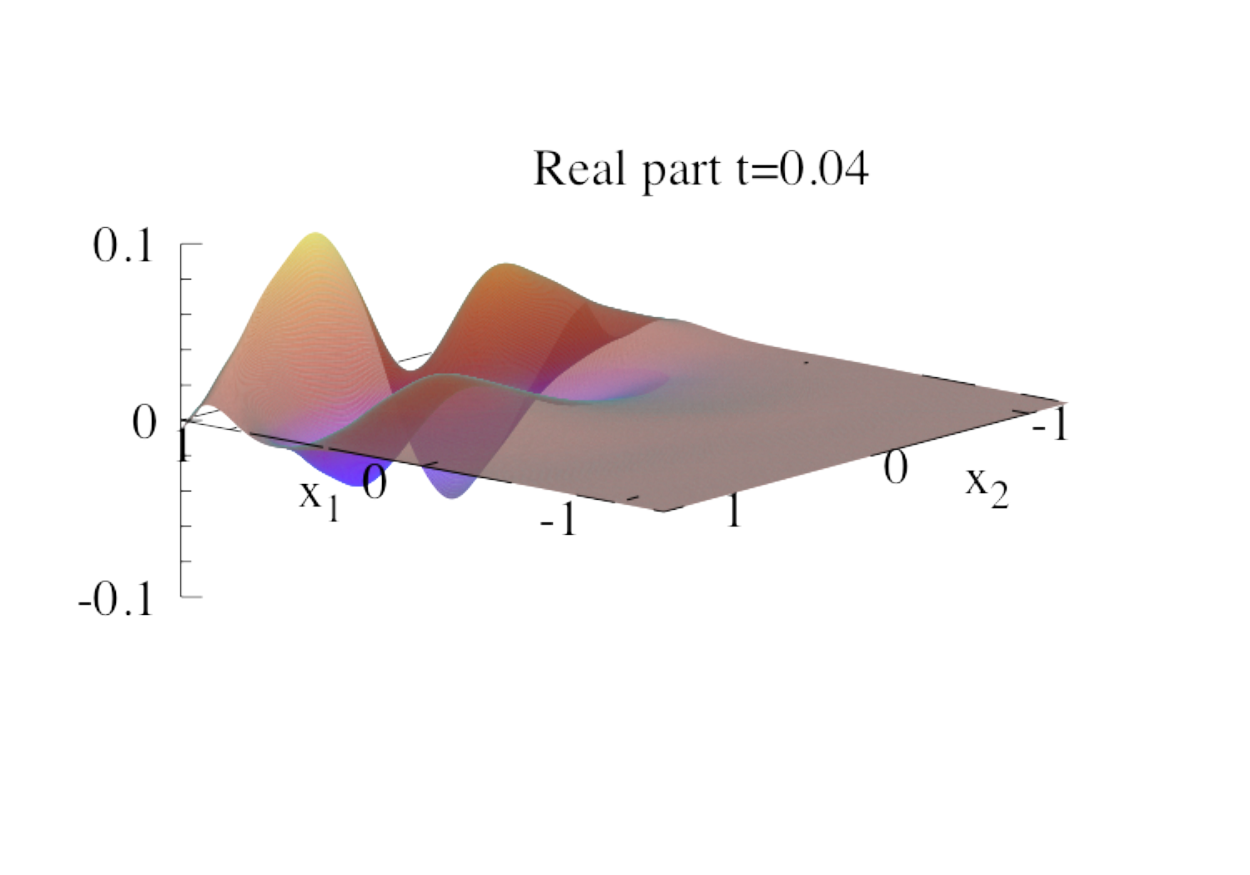} 
\vspace{-1.2cm}
\caption{Real part of $u$ when $w_1(x)=\ee^{30i x} \ee^{-60(x-1/4)^2}$, $w_2(x)=\sin(\pi x)$,  $N=6$, $h=0.005$\,.}\label{figresin}
\includegraphics[scale=0.5]{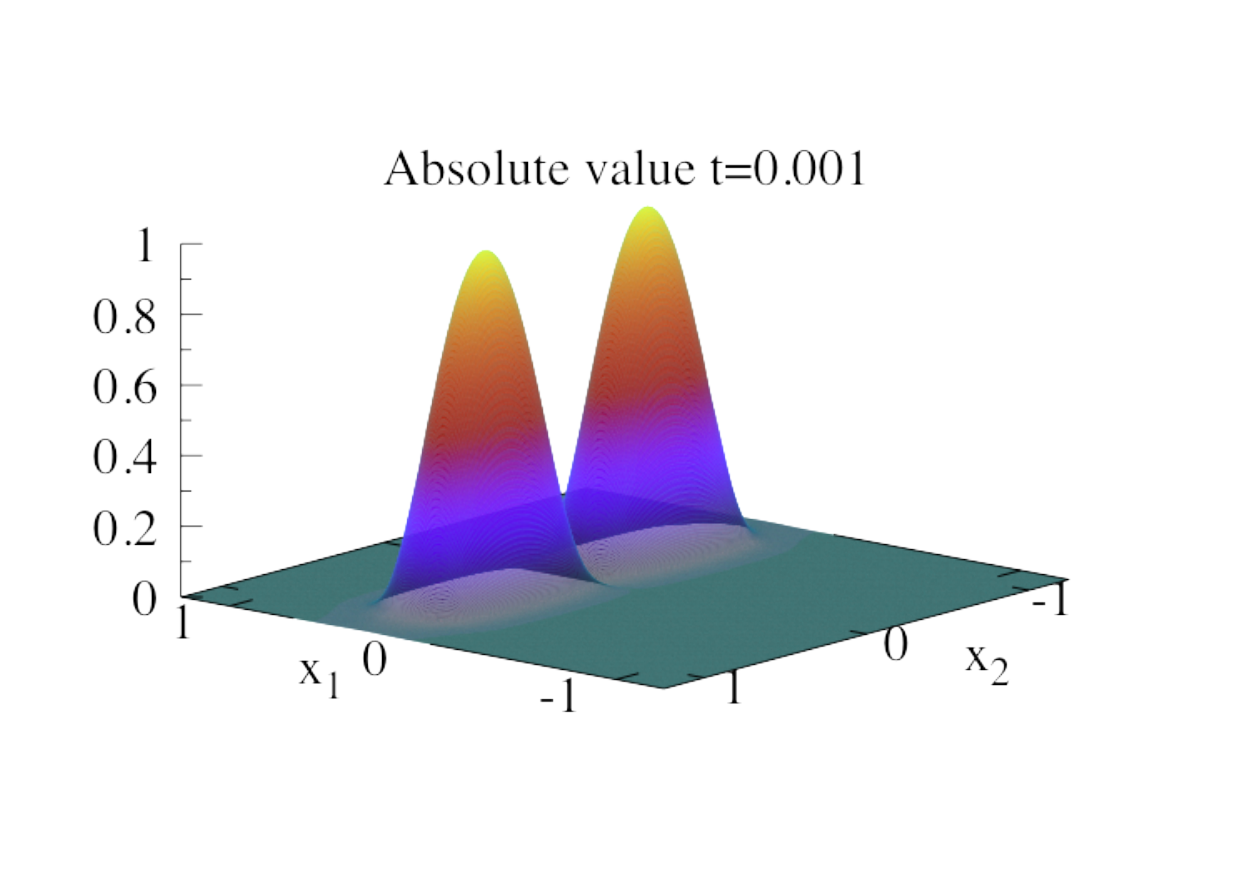} 
\includegraphics[scale=0.5]{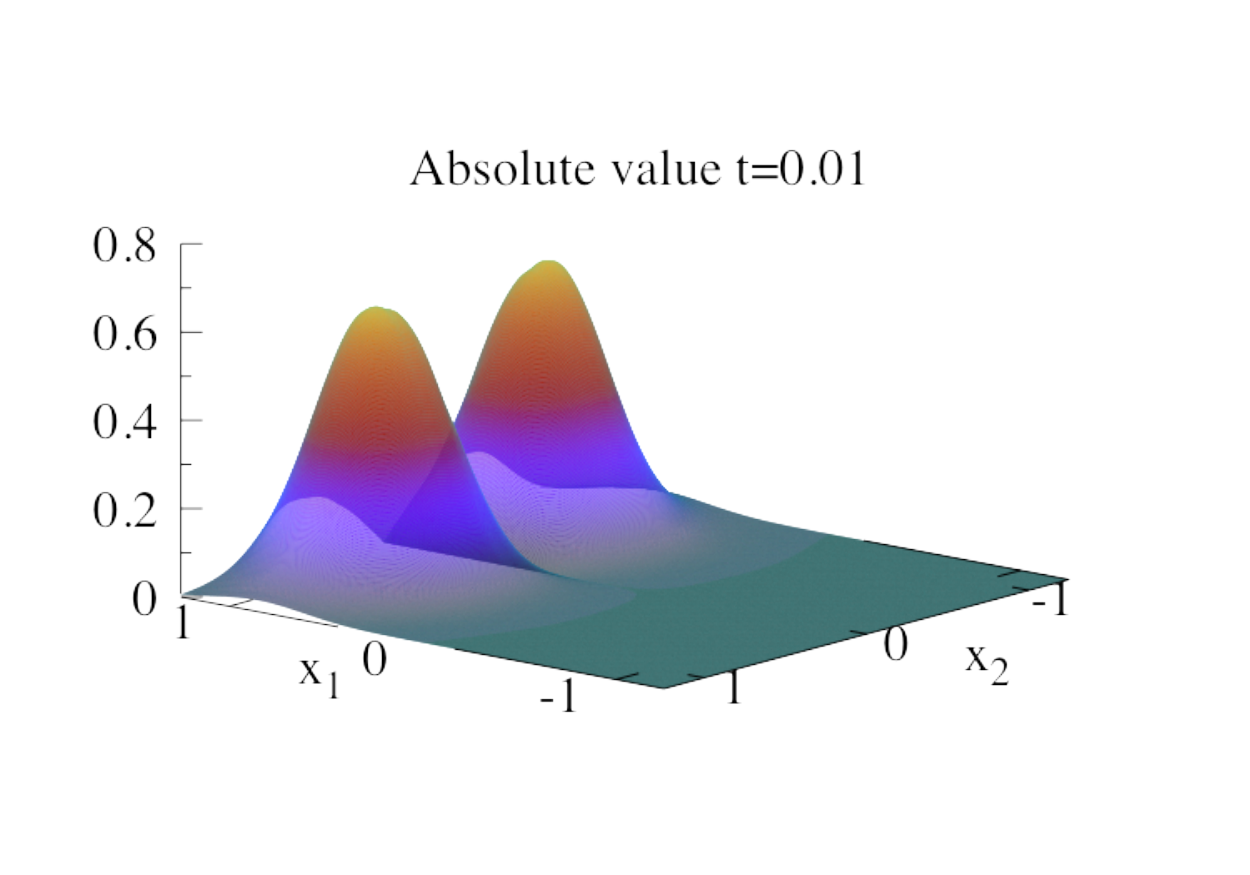}\\[-1cm]
\includegraphics[scale=0.5]{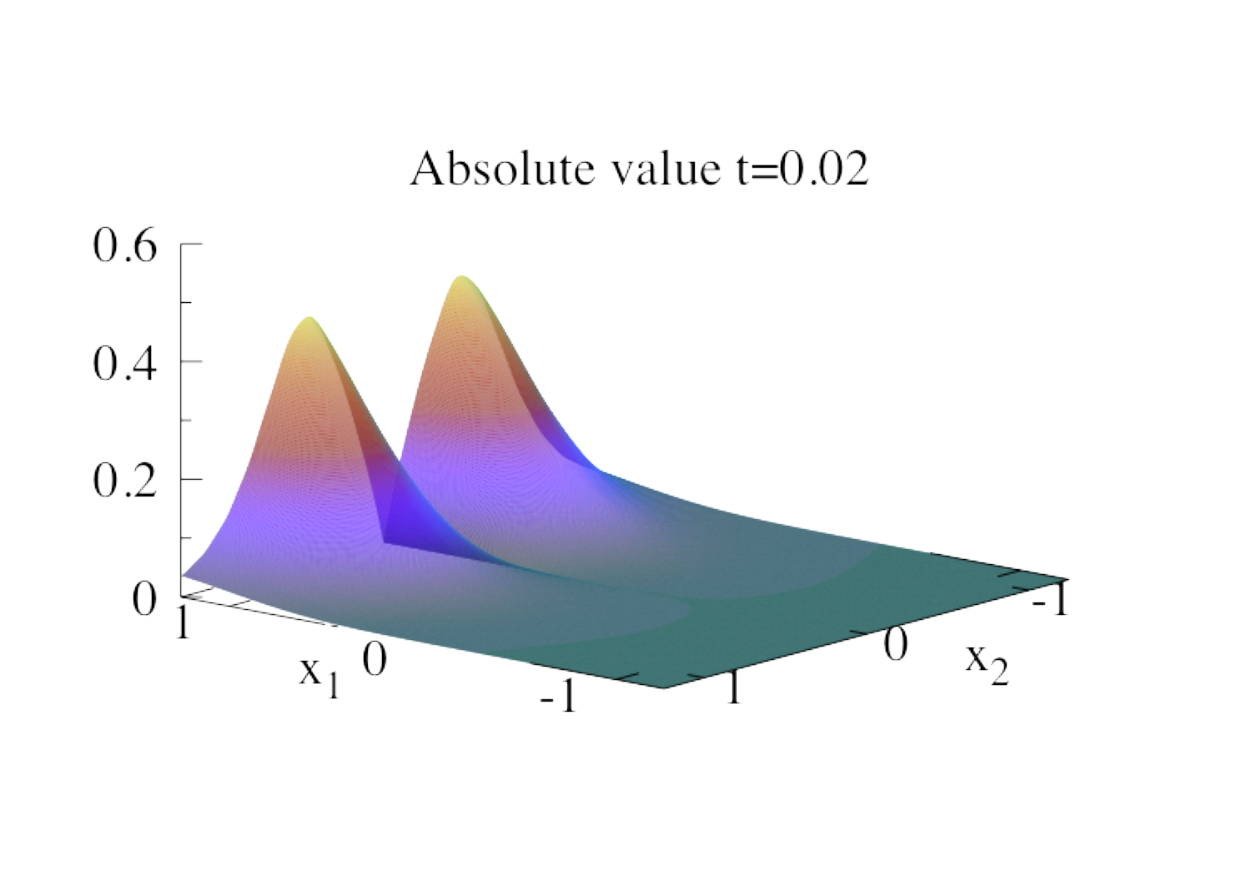}  
\includegraphics[scale=0.5]{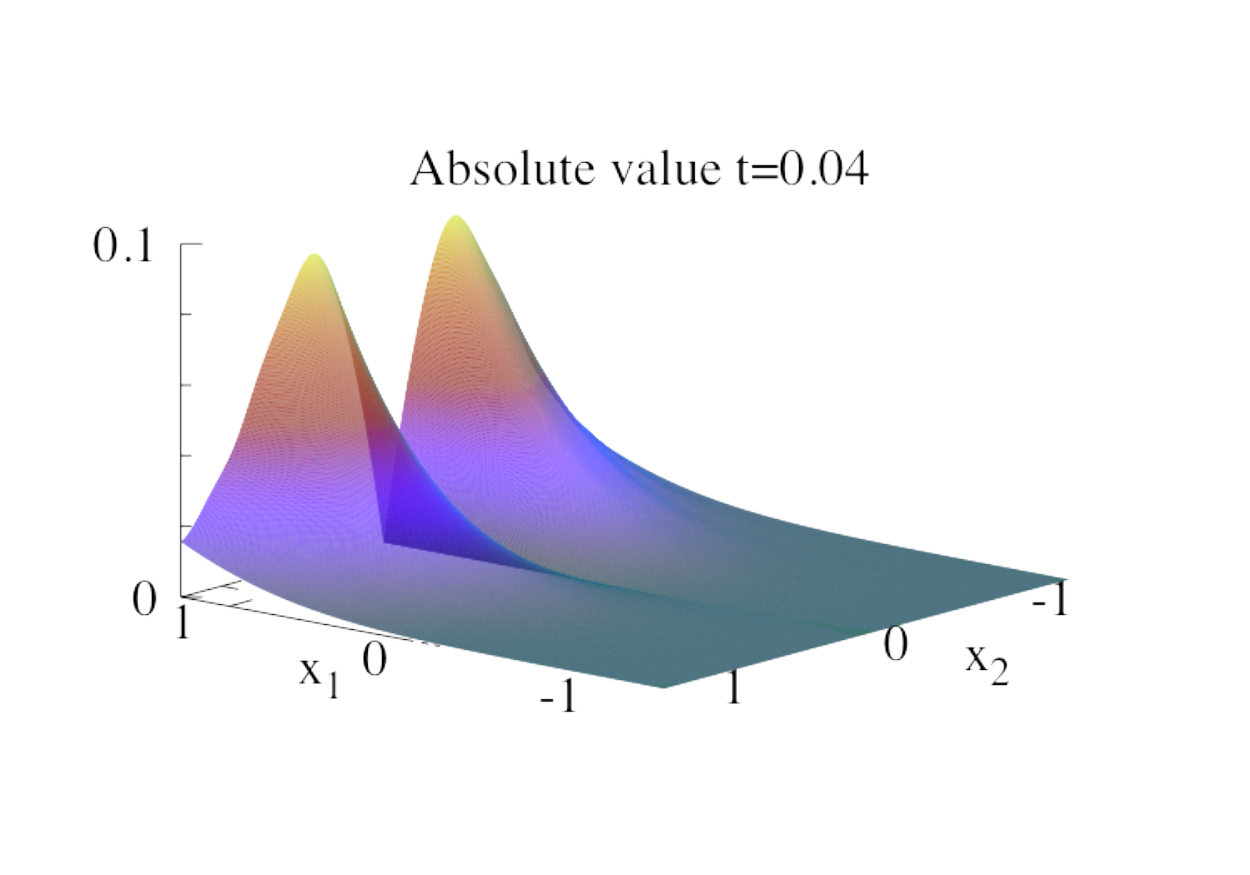} 
\vspace{-1.2cm}
\caption{Absolute value of $u$ when $w_1(x)=\ee^{30i x} \ee^{-60(x-1/4)^2}$, $w_2(x)=\sin(\pi x)$,  $N=6$, $h=0.005$\,.}\label{figabsin}
\end{figure}

We conclude the paper illustrating the evolution of  $u(\bx,t)$
evolving under the two-dimensional Schr\"odinger equation
\eqref{ex1}. 
 First we consider the evolvement of the traveling
Gaussian
$
u(\bx,0)=\ee^{c\, i (x_1-x_2)} \ee^{-60|\bx|^2} 
$
on the domain $(-1.25,1.25)\times (-1.25,1.25)$ at four consecutive
time values. 
Figures \ref{figre30} and \ref{figab30} show the %temporary
evolution of $\re u(\bx,t)$ and $|u(\bx,t)|$ when $c=30$. At time
$t=0.04$ the solution has almost completely left the domain. The case
$c=10$ is reported in Figures \ref{figre10} and \ref{figab10}.
Figures \ref{figresin} and \ref{figabsin} concern the initial data
$
g(x_1,x_2)=\ee^{30\,i  x_1} \ee^{-60(x_1-1/4)^2}\sin(\pi x_2)\,.
$
The figures of the imaginary part of $u(\bx,t)$ are
virtually the same as for the real part, so we skipped that plots. In
all the figures we used the approximation formula of order $N=6$, the
extension of $\widetilde{w}_j=w_j$ and the step $h=0.005$. 
Similar tests with finite difference scheme can be found in \cite{AS2} and \cite{SA1}.

\end{document}